\documentclass [12pt]{amsart}

\usepackage{amsmath, accents} %

\usepackage{amssymb,amsxtra,amsfonts}
\usepackage{epsf}              %
\usepackage{epsfig}             %
\usepackage{graphics}
\usepackage{epigraph}

\usepackage{tikz}
\usepackage{wrapfig}

\usepackage[colorlinks]{hyperref}
\usepackage{combelow}   %

\usepackage{scalerel}   %
\def\apple{{\includegraphics[height=2ex]{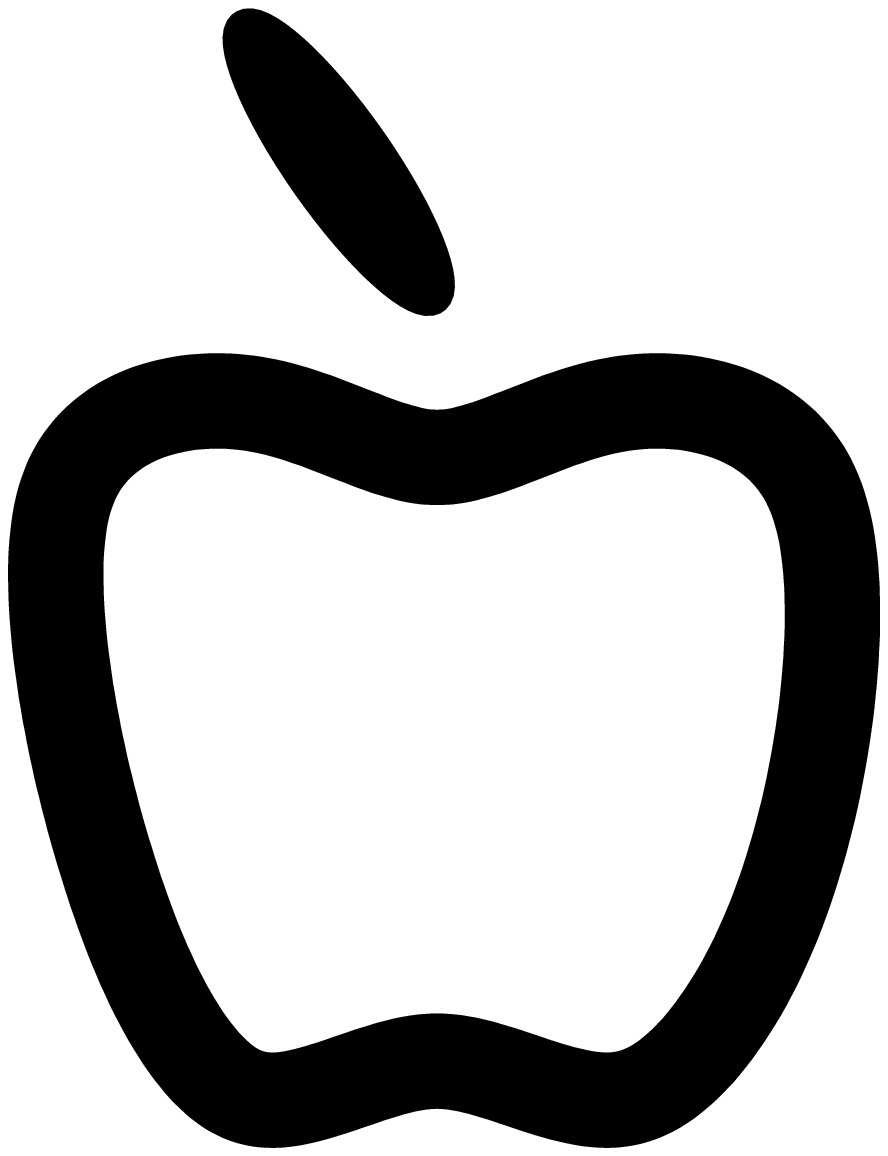}}}

\usepackage{etoolbox}

\makeatletter
\patchcmd{\@setauthors}{\MakeUppercase}{}{}{}
\makeatother

\long\def\pb #1*/{}

\def\reE@DeclareMathSymbol#1#2#3#4{%
    \let#1=\undefined
    \DeclareMathSymbol{#1}{#2}{#3}{#4}}
\DeclareSymbolFont{symbolsC}{U}{txsyc}{m}{n}
\SetSymbolFont{symbolsC}{bold}{U}{txsyc}{bx}{n}
\DeclareFontSubstitution{U}{txsyc}{m}{n}
\reE@DeclareMathSymbol{\strictiff}{\mathrel}{symbolsC}{76}

\newcommand\beq{\begin{equation}}
\newcommand\eeq{\end{equation}}
\newcommand\bal{\begin{align*}}
\newcommand\eal{\end{align*}}   
\newcommand\bmx{\left(\begin{matrix}}
\newcommand\emx{\end{matrix}\right)}
\newcommand\bsmx{\left(\begin{smallmatrix}}
\newcommand\esmx{\end{smallmatrix}\right)}
\newcommand\bmxnp{\begin{matrix}}
\newcommand\emxnp{\end{matrix}}
\newcommand\bsmxnp{\begin{smallmatrix}}
\newcommand\esmxnp{\end{smallmatrix}}

\newcommand{\wt}{\widetilde}

\DeclareMathSymbol{\widehatsym}{\mathord}{largesymbols}{"62}
\newcommand\lowerwidehatsym{%
  \text{\smash{\raisebox{-1.3ex}{%
    $\widehatsym$}}}}
\newcommand\fixwidehat[1]{%
  \mathchoice
    {\accentset{\displaystyle\lowerwidehatsym}{#1}}
    {\accentset{\textstyle\lowerwidehatsym}{#1}}
    {\accentset{\scriptstyle\lowerwidehatsym}{#1}}
    {\accentset{\scriptscriptstyle\lowerwidehatsym}{#1}}
}

\newcommand{\wh}{\fixwidehat}

\newcommand{\bSi}{{\bf \Si}}

\newcommand{\monesp}{\;\!\!}   
\newcommand{\onto}{\twoheadrightarrow}
\newcommand{\into}{\hookrightarrow}

\newcommand{\st}{\ \bigl\vert\ }

\providecommand{\Rad}{\text{\rm Rad}}

\providecommand{\from}{\leftarrow}

\providecommand{\<}{\langle}
\renewcommand{\>}{\rangle}

\def\part#1{\frac{\partial\phantom{q}}{\partial#1}}

\newcommand {\flp}{(\!(}
\newcommand {\frp}{)\!)}

\newcommand{\sdp}{{\rtimes}}    

\newcommand{\HH}{\text{\rm H}}

\newcommand{\Lie}{{\mathop{\rm Lie}}}

\DeclareMathOperator{\Sto}{{\cS}to} 
 
\DeclareMathOperator{\ISto}{{\IS}to} 
 
\newcommand{\Sect}{\text{\rm Sect}}

\DeclareMathOperator{\rk}{\mathop{\rm rk}}
\newcommand{\ram}{\mathop{\rm Ram}}

\newcommand{\Splits}{{\rm Splits}}

\DeclareMathOperator{\Iso}{Iso}   
\DeclareMathOperator{\GrIso}{GrIso}   
\DeclareMathOperator{\GrAut}{GrAut}

\DeclareMathOperator{\Hom}{Hom}         
\DeclareMathOperator{\Aut}{\mathop{\rm Aut}}

\newcommand{\GL}{{\mathop{\rm GL}}}

\newcommand{\Gr}{{\rm Gr}}

\DeclareMathOperator{\Rep}{\rm Rep}

\renewcommand{\ker}{\mathop{\rm Ker}}

\DeclareMathOperator{\End}{End}

\newcommand{\diag}{{\mathop{\rm diag}}}

\newcommand{\hk}{{hyperk\"ahler }}   

\newcommand{\ba}{{\bf a}}
\newcommand{\bA}{{\bf A}}

\newcommand{\bB}{{\bf B}}

\newcommand{\bd}{{\bf d}}

\newcommand{\bH}{{\bf H}}

\newcommand{\bS}{{\bf S}}

\DeclareSymbolFont{bbold}{U}{bbold}{m}{n}
\DeclareSymbolFontAlphabet{\mathbbold}{bbold}

\newcommand{\IA}{\mathbb{A}}

\newcommand{\IC}{\mathbb{C}}

\newcommand{\IF}{\mathbb{F}}

\newcommand{\IH}{\mathbb{H}}

\newcommand{\IN}{\mathbb{N}}

\newcommand{\IP}{\mathbb{P}}                                     
\newcommand{\IQ}{\mathbb{Q}}                           
\newcommand{\IR}{\mathbb{R}}                           
\newcommand{\IS}{\mathbb{S}}

\newcommand{\IT}{\mathbb{T}}

\newcommand{\IV}{\mathbb{V}}

\newcommand{\IZ}{\mathbb{Z}}

\newcommand{\cB}{\mathcal{B}}

\newcommand{\LocSys}{{\mathcal{L}oc\,\monesp\monesp\mathcal{S}ys}}

\newcommand{\SLocSys}{{\mathcal{S}\monesp\monesp\mathcal{L}oc\,\monesp\monesp\mathcal{S}ys}}

\newcommand{\cF}{\mathcal{F}}
\newcommand{\cG}{\mathcal{G}}

\newcommand{\cI}{\mathcal{I}}

\newcommand{\cM}{\mathcal{M}}

\newcommand{\cR}{\mathcal{R}}
\newcommand{\cS}{\mathcal{S}}

\newcommand{\cT}{\mathcal{T}}

\newcommand{\cV}{\mathcal{V}}

\newcommand{\gM}{       \mathfrak{M}     }

\newcommand{\gl}{       \mathfrak{gl}     } 

\newcommand{\al}{\alpha}

\newcommand{\be}{\beta}
\newcommand{\ga}{\gamma}

\newcommand{\De}{\Delta}

\newcommand{\Ga}{\Gamma}

\newcommand{\la}{\lambda}
\newcommand{\La}{\Lambda}

\newcommand{\Om}{\Omega}
\newcommand{\Si}{\Sigma}
\newcommand{\Th}{\Theta}
\renewcommand{\th}{\theta}

\newcommand{\ze}{\zeta}

\renewcommand{\bar}{\overline}

\makeatletter
 \newlength{\typesize}
\makeatother

\newlength{\vvoff}
\newlength{\hhoff}

\def\mapright#1{\smash{
        \mathop{\longrightarrow}\limits^{#1}}}

\newcommand{\pf}{\begin{bpf}}

\newcommand{\pfms}{\begin{bpfms}}
\newcommand{\epf}{\end{bpf}\hfill$\square$\\}           
\newcommand{\epfms}{\end{bpfms}\hfill$\square$\\}       

\newcommand{\idea}{\begin{bidea}}

\newcommand{\eidea}{\end{bidea}\hfill$\square$\\}           

\newcommand{\sk}{\begin{bsk}}    

\newcommand{\esk}{\end{bsk}\hfill$\square$\\}           
\newcommand{\sketch}{\begin{bsketch}}

\newcommand{\esketch}{\end{bsketch}\hfill$\square$\\}

\newtheorem {hypo}{\bf\hspace{-\parindent}Hypothesis}
\newtheorem {thm}[hypo]{Theorem}   
\newtheorem {prop}[hypo]{Proposition}

\newtheorem {cor}[hypo]{Corollary}
\newtheorem {lem}[hypo]{Lemma}

\theoremstyle{definition}\newtheorem {defn}[hypo]{Definition}
\theoremstyle{definition}
\theoremstyle{definition}\newtheorem{eg}[hypo]{Example} 
\theoremstyle{remark}\newtheorem{rmk}[hypo]{Remark}
\theoremstyle{remark}

\numberwithin{equation}{section}  %
\numberwithin{hypo}{section}

\DeclareMathAlphabet{\mathbbmsl}{U}{bbm}{m}{sl}

\begin{document}

\pagenumbering{gobble}

\title[Topology of the Stokes phenomenon]{Topology of the Stokes phenomenon}
\author[Philip Boalch]{P. P. BOALCH\\ \  \\
Institut de Math\'ematiques de Jussieu -- Paris Rive Gauche\\
Universit\'e de Paris \& CNRS \\ 
B\^atiment Sophie Germain \\ 
8 Place Aur\'elie Nemours \\
75205 Paris, FRANCE}

\begin{abstract}
Several intrinsic topological ways to encode 
connections on vector bundles on 
smooth complex algebraic curves will be described.
In particular the notion of {\em Stokes decompositions}
will be formalised,  
as a convenient intermediate category between 
the Stokes filtrations and the Stokes local 
systems/wild monodromy representations.
The main result establishes a new simple characterisation of the Stokes decompositions.

\ 

\  

\

$$\ \ \ \includegraphics[width=0.6\textwidth]{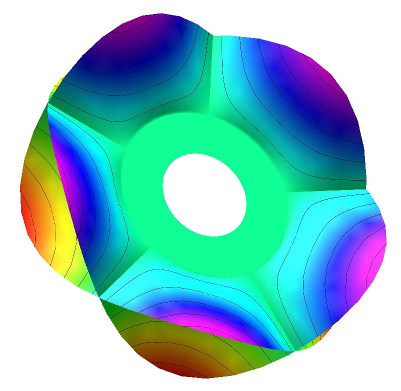}
$$

\end{abstract}


\maketitle

\newpage

\pagenumbering{arabic}
\setcounter{page}{2}

\ 

\

\epigraph{
``The subject ought to be one of pure mathematics, for
it is in honour of ABEL, and most of my work refers to applications of mathematics. There is one thing I thought might perhaps do....  The subject is the discontinuity of arbitrary constants that appear as multipliers of semi-convergent series...''

\ 

G.G. Stokes, Cambridge, 23/4/1902
}

\setcounter{tocdepth}{1}
\tableofcontents

\section{Introduction}

\subsection{} 
Systems of meromorphic linear differential equations on the complex plane have been studied for centuries, and, more generally, algebraic connections on vector bundles on smooth complex algebraic curves are extremely basic objects.

The {\em Riemann--Hilbert correspondence} 
\cite{Del70, Katz76} says that the special class of connections with regular singularities (moderate growth/Fuchsian type) are classified by their local system of solutions, i.e. by their monodromy representations (upon choosing a basepoint).
If one sets this up carefully (as in \cite{Del70})
it is an equivalence of categories---this might sound exotic but it really is just a very convenient way to express precisely the fact that the topological data (the local system/fundamental group representation) encodes the algebraic differential equations.
In particular it implies the more familiar fact that two connections are isomorphic if and only if their monodromy representations are conjugate.

For more general connections, those with {\em irregular singularities}, the situation is unfortunately much less widely understood in the mathematics community.
The basic question is to describe the data that one needs to add to the local system to encode such connections.
The fact of the matter is that this can indeed be done, and there are several different approaches, all of which are equivalent, but not in a trivial way.
A basic difficulty is to recognise the intrinsic geometric structures that appear---one needs to be able to see through the explicit coordinate dependent computations prevalent (and essential) in the development of the subject.
In the regular singular case this was done a long time ago and underlies the definition of the fundamental group and the notion of local system (locally constant sheaf of complex vector spaces).

The aim of this article is to describe intrinsically, as simply as possible, 
three different approaches to the extra ``topology at infinity''
that arises in the case of irregular connections:
1) {\em Stokes filtrations}, 2) {\em Stokes decompositions/gradings}, and 3) {\em Stokes local systems}, or {\em wild monodromy representations}.
We will then give a new direct proof (building on an idea of Malgrange) that they are equivalent.
The reader will see that the situation is somewhat analogous to Hodge structures, which may be described equivalently in terms of the Hodge filtration or the Hodge decomposition.

The impetus came from studying a simple class of examples (in \cite{even-euler}) and then trying 
understand directly the equivalence of 1) and 3) in general: 
it turns out that this becomes much simpler 
if one first formalises the category 2) of Stokes graded local systems, and then 
establishes two equivalences 
1) $\xleftarrow{\raisebox{-0.7ex}[0ex][0ex]{$\,\sim\,$}}$ 
2) $\xrightarrow{\raisebox{-0.7ex}[0ex][0ex]{$\,\sim\,$}}$ 3), 
as will be done here.

The rest of this introduction will give a sketch, deferring full definitions to later sections.
Some of the recent motivating questions 
(construction of wild character varieties, 
the nonlinear local systems that they form, the relative Riemann--Hilbert--Birkhoff maps, the link to Drinfeld--Jimbo quantum groups, 
and the TQFT approach to meromorphic connections) are reviewed in 
\cite{smid, HDR, hit70, ISTtalks}. 
In particular the wild character varieties 
(moduli spaces of Stokes local systems) 
give the simplest description of the differentiable manifolds underlying a somewhat vast 
collection of complete \hk manifolds, 
the {nonabelian Hodge spaces
(see \cite{hit70} Defn. 7).
In a 
different complex structure they are algebraically completely integrable Hamiltonian systems,
some key examples of which arise as finite gap solutions of integrable hierarchies, such as KdV \cite{DMN76}.
On the other hand, as noted in \cite{smid},
an initial impetus for these questions 
came from Dubrovin's classification of
semisimple Frobenius manifolds in terms of Stokes data.

\subsection{} 
One completely general  approach 
involves adding filtrations (flags indexed by certain ordered sets) on sectors at each pole, the {\em Stokes filtrations}.
The exact types of filtrations that occur, and how they may jump,
can be axiomatised yielding an equivalence of categories 
\cite{deligne78in, Mal-irregENSlong}, generalising the Riemann--Hilbert correspondence.
The filtrations encode the exponential growth rates of solutions as one goes towards a pole.
(The approach in \S\ref{sn: SFLS} below is very close to 
op. cit., 
but slightly simpler since  we index the filtrations by finite totally ordered sets.)
A consequence of this is that there is a canonically defined (continuous) {\em graded local system} on a small punctured disk around each pole: the associated graded of the Stokes filtrations.
This is an intrinsic way to describe the formal classification (of Fabry--Hukuhara--Turrittin--Levelt): for example two connections germs are formally isomorphic (over $\IC\flp z\frp$) if and only if their 
graded local systems are isomorphic.

The revolutionary idea of Stokes \cite{stokes1857} was that it is sometimes better to start at infinity, at the pole, and work from there towards the interior of the curve.
He used this idea to give a spectacular computation 
of the position of the fringes of rainbows (the zeros of the Airy function) far more efficiently than summing the Taylor series at zero.
Mathematically a key idea appearing here is that: 

a) the connection canonically determines a finite number of ``singular directions'' at each pole, and 

b) the choice of any formal solution of the connection
at a pole uniquely determines a preferred actual 
solution, on any sector at the pole not containing a singular direction.

Stokes only worked this out in a few examples (related to the Airy and Bessel equations), but the idea is correct and is now a general theorem (multisummation of formal solutions of linear differential equations \cite{BBRS91}).
It also holds more generally in many non-linear situations when one does a formal simplification
\cite{Ecalle81, MR91}.

The import of this to the problem of describing the category of connections topologically, 
is the following:
\beq \label{eq: stokes gluing}
\begin{matrix}
\text{\em Away from the singular  directions there is a preferred/canonical way to glue}\\
\text{\em the graded local system to the local system of solutions of the connection.}
\end{matrix}
\eeq

In particular in these sectors there is a preferred grading, the {\em Stokes grading} or {\em Stokes decomposition} (not just a filtration) of the local system of solutions.
One can axiomatise the gradings that occur and 
their possible discontinuities across the singular 
directions, 
and again prove 
an equivalence,  between connections and Stokes graded local systems.

The Stokes gradings split the Stokes filtrations wherever they are both defined, and this actually characterises them:

\begin{thm}[Splitting] \label{thm: mainresult}
For any Stokes filtered local system there is a unique 
Stokes graded local system (with the same underlying local system)
such that the gradings split the filtrations wherever both are defined. 
\end{thm}

This statement looks to be new. See \S\ref{ssn: summary} and Thm.  \ref{thm: v2 of main thm}  for more details, and 
Prop. \ref{prop: one level Stokes grading} for the simplest case, with just one level. 
Note that usually one gets preferred gradings
by having fixed asymptotics (i.e. splitting the Stokes filtrations) on a large enough sector 
(or a nested sequence of such sectors), 
whereas the characterisation here is different.
On the other hand one can axiomatise the local systems that occur when the graded local system
is actually glued  
to the local system of solutions. 
One way to formalise this is to boldly puncture 
the underlying curve $\Si^\circ$ near each singular direction at each singular point, to yield a new curve 
$\wt \Si\subset \Si^\circ$. 
The idea underlying these {\em tangential punctures} is already in Stokes' paper---see 
\S\ref{sn: borel sums} for more on this.
The Stokes gluing \eqref{eq: stokes gluing} then
yields a local system on $\wt \Si$, the {\em Stokes local system} \cite{gbs, twcv}, equal to the graded local system on a halo (small punctured disk) near each pole, and to the
local system of solutions away from the halos.
The Stokes local systems can be axiomatised and again one
gets an equivalence of categories (it is very close to the category of Stokes graded local systems).

The resulting equivalence between Stokes filtrations and Stokes local systems could be viewed as an intrinsic global version
of the main theorem of Loday-Richaud \cite{L-R94} 
(no longer needing the choice of a 
marking/formal normal form, 
or any discussion of nonabelian cocycles) 
which in turn refines the main result 
(Thm. 3.4.1) of Babbitt--Varadarajan \cite{BV89}
and is essentially equivalent to results of Jurkat in 
\cite{jurkat78, BJLproper}.

Thus in summary this gives three topological descriptions of the category of connections: 
as Stokes filtered local systems, Stokes graded local systems or as Stokes local systems. 
Combined, they give a multi-faceted answer to the question ``What is a Stokes structure?'', much as a Hodge structure has two different descriptions.

Whereas the Deligne--Malgrange approach to Stokes filtrations is easily seen to be intrinsic (once one realises the exponential local system has an intrinsic definition, as in \cite{twcv} Rmk 3), the other approaches are less well developed,  
so we are taking the opportunity here
to describe them intrinsically, 
thus giving a topological approach to the version of the Stokes phenomenon actually discovered by Stokes (hence the title of this article).

In turn these three approaches give 
three different approaches to the wild character 
varieties: just as moduli spaces of local systems form  an interesting class of varieties (the character varieties), moduli spaces of Stokes local systems form the wild character varieties.
Due to the three descriptions, each wild character variety
solves three different moduli problems: it is also a moduli space of Stokes filtered local systems, and 
a moduli space of Stokes graded local systems 
(cf. e.g. the examples in \cite{even-euler, sibuya1975}).

This is very useful in practice since the description in terms of Stokes local systems yields an explicit presentation of the wild character variety as a quotient 
\beq\label{eq: WCV}
\Hom_\IS(\Pi,G)/\bH,\qquad \Pi=\Pi_1(\wt \Si,\be)
\eeq
of an affine variety $\Hom_\IS(\Pi,G)$ (the space of Stokes representations, or wild monodromy representations), by a reductive group $\bH$.
This enables the use of standard geometric invariant theory to construct the wild character variety
algebraically \cite{saqh, gbs, twcv}.
The possibility of such a presentation is not evident if one starts from the Stokes filtration viewpoint.
Further, earlier approaches use a different type of framing called a ``marking'', and in general this does not lead to reductive groups---see \cite{BV89} Theorem 2.3.1, or \cite{BV85-bams} p.98.
(So in essence we have rejigged things in order
to be able to  carry out this construction.)
This is a direct generalisation of the standard presentation of the 
(tame) character variety in the form
\beq\label{eq: CV}
\Hom(\pi_1(\Si^\circ,b),G)/G.
\eeq
A major preoccupation has been to 
generalise geometric properties of the spaces 
\eqref{eq: CV} to the spaces \eqref{eq: WCV}, such as:  1) symplectic/Poisson structures 
\cite{FN82, Ugag, thesis, smid, saqh, gbs, twcv}
(the spaces \eqref{eq: WCV}
have algebraic Poisson structures with symplectic leaves given by fixing the isomorphism classes of the graded local systems at each pole), %
2) hyperk\"ahler metrics and special Lagrangian fibrations 
\cite{wnabh, Sab99}
(sufficiently generic symplectic leaves are complete hyperk\"ahler manifolds, becoming meromorphic Hitchin systems in special complex structures), and
3) wild mapping class group actions
\cite{Mal-imd12, bafi, gbs}
(upon deforming the underlying wild Riemann surface
the spaces \eqref{eq: WCV} form a local system of Poisson varieties, and the monodromy of this is the
wild mapping class group action on \eqref{eq: WCV}).

\subsection{Summary of main result}\label{ssn: summary}
In brief the logic of Thm. \ref{thm: mainresult} is as follows.

1) The Stokes graded local systems involve adding 
gradings to a local system
on sectors at each puncture.
This involves data $n,I,\Th,\IA,\prec_d$
where:

\noindent$\bullet$\ 
$n$ is the rank of the local systems,

\noindent$\bullet$\ 
$I$ is a finite covering of the circle of directions at each pole (used to index the gradings),

\noindent$\bullet$\ 
$\Th$ an integer for each component of $I$
(the dimensions of the graded pieces),

\noindent$\bullet$\ 
$\IA$ is a finite set of directions at each puncture, where the gradings may jump,

\noindent$\bullet$\ 
$\prec_d$ (the Stokes arrows) 
is a partial order of $I_d$ for each 
$d\in \IA$ that will be used to control
the possible relative positions of the gradings across each direction in $\IA$.

2) Similarly the Stokes filtered local systems involve adding 
filtrations
on sectors at each puncture.
This involves data $n,I,\Th,\IS,<_d$
where $n,I,\Th$ are as above and:

\noindent$\bullet$\ 
$\IS$ is a finite set of directions at each puncture, where the filtrations 
may jump,

\noindent$\bullet$\ 
$<_d$ (the exponential dominance orderings) 
is a total order of $I_d$ for each direction
$d\not\in \IS$ (continuous provided $d$ does not cross $\IS$). They will be used to index the filtrations and to control
the possible relative positions of the filtrations across each direction in $\IS$.

3) The notion of irregular class will be recalled 
in \S\ref{sn: irclass}. The key point is that it canonically determines both types of data: 
$(n,I,\Th,\IA,\prec_d)$ 
and $(n,I,\Th,\IS,<_d)$.
(Also any algebraic 
connection on the curve canonically determines an irregular class at each puncture.)
Thus one can consider Stokes graded local systems
and Stokes filtered local systems with the given irregular class.
The main result then says that for each 
Stokes filtered local system 
of type $(n,I,\Th, \IS,<_d)$, 
there is a {\em unique}
Stokes grading of type 
$(n,I,\Th, \IA,\prec_d)$ (on the same underlying local system) that splits the Stokes filtration wherever both are defined.

This result is simple in the rank two cases studied classically (second order equations): the Stokes filtrations are then given by a line (the subdominant solutions) 
in a rank two local system, and
across $\IS$ such lines 
are transverse and make up the Stokes grading.
See \cite{even-euler} for a recent exposition---as 
shown there even 
in such  examples the resulting different descriptions of the wild character variety are interesting, for example relating the Euler continuant polynomials to the multiratio of tuples of points on the Riemann sphere.
The simplest of the examples in \cite{even-euler} 
will be used as a running example  
to illustrate some of the basic constructions, 
although it doesn't exhibit the complexity of the general set-up.

\subsection{Further background/other approaches}\label{ssn: further background}

The history is long and complicated and would require a book to put it all in its proper context. 
Note that most authors work entirely in one point of view. 
Nonetheless here is a brief attempt 
to better document the 
various discoveries leading to the  wonderful fact that
the category of connections on smooth curves has 
a precise topological description (or more precisely, several). Cf. also \cite{ramis-history, Var96, DMRci} and references therein.

1) Similar presentations to \eqref{eq: WCV}
were first 
found by Birkhoff \cite{birkhoff-1909,birkhoff-1913}
for generic connections of any rank.
This work 
really started the subject of constructing 
invariants of irregular connections.
It was extended from the generic case to 
the general case 
in \cite{jurkat78, BJLproper}, using the 1976 version of 
\cite{SibBk}.
Their approach also involves preferred bases on sectors. 
It is similar but not quite the same as the Stokes approach used here.
An intrinsic form of their results 
will be discussed in a sequel to the present paper.
Whereas Stokes filtrations arise as one goes 
{\em towards} a pole, 
and Stokes gradings arise if one starts {\em at} a pole, Birkhoff analysed what happens when one goes {\em around} a pole.

2) The Malgrange--Sibuya cohomological approach \cite{Mal79, sibuya77} (cf. \cite{BV89, L-R94, SibBk}) 
looks to have been
the first general classification 
of germs of marked connections, and so 
it traditionally acts as a 
hub to pass between different viewpoints
(or as solid ground on which to establish other viewpoints). 
The term {\em Stokes structure} first appeared in 
\cite{sibuya77}.
As emphasised above, the message of 
\cite{deligne78in, Mal-irregENSlong} is that it is sometimes better to think in terms of filtrations 
(early examples of which are the subdominant solutions in \cite{sibuya1975}), 
and  the message of \cite{stokes1857} 
(and \cite{birkhoff-1909, MR91, L-R94})
is that it is sometimes better to think in terms of gradings, or  wild monodromy.

3) The Stokes approach goes 
back to Stokes' paper \cite{stokes1857}. 
Its extension to generic connections of arbitrary rank is essentially the story explained 
in \cite{BJL79}
(this is similar to 
\cite{birkhoff-1909,birkhoff-1913}, but solved the 
{\em Birkhoff wall-crossing problem}, of central importance in isomonodromy).
It appears in many works on isomonodromy and the Painlev\'e equations
such as \cite{JMU81, smid, fikn}.
The general case appears (in slightly different forms) in Martinet--Ramis \cite{MR91} and Loday-Richaud \cite{L-R94} (and the links to 
work of Ecalle and others on multisummation are explained 
there---pioneering work of E.Borel, Watson and Dingle underlies this approach). 
These provided much inspiration 
and the main result above (Thm. \ref{thm: mainresult}) 
resulted from trying to understand intrinsically 
the main theorem of 
\cite{L-R94} (describing the Malgrange--Sibuya non-abelian cohomology space in terms of Stokes groups).

\subsection{Generalisations}
Some extensions to be discussed in detail elsewhere are as follows.

1) In modern applications of these results, such as 
isomonodromy/wall-crossing or 2d gauge theory (wild non-abelian Hodge theory on curves), a
slightly bigger category than the category of connections on curves is used, involving particular extensions across the punctures 
(and compatible parabolic structures/filtrations, cf. 
\cite{wnabh, ihptalk, hit70}). Topologically this involves upgrading each graded piece of the graded local system to be  $\IR$-filtered.
This was understood in the tame case in 
\cite{Sim-hboncc} ($\IR$-filtered local systems) and is not essentially different in the wild case ($\IR$-filtered Stokes local systems),
once one understands how to 
superpose  the tame story 
on the wild story.
The bijective Riemann--Hilbert--Birkhoff correspondence of \cite{smid} Cor. 4.9 %
is an example of this.

2) %
The work 
\cite{bafi, saqh, fission, gbs, twcv} involves the extension of 
Stokes data to the case of connections on 
principal $G$-bundles
for other algebraic groups beyond $\GL_n(\IC)$,
mainly from the viewpoint of Stokes local systems.
It requires some familiarity with root systems/Lie theory.
The original aim was to write the present article 
at that level of generality
but the simple uniqueness statement of Thm. \ref{thm: mainresult}  did not seem to be known even for 
$\GL_n(\IC)$, so a separate (simpler) presentation of this case seemed justified.

3) Of course having a clear picture of how the linear case works, suggests how to phrase the {\em nonlinear} Stokes phenomenon 
(this is already remarked in \cite{MR91}).
More pointedly the 
version for general principal $G$-bundles suggests directly how the nonlinear case should work,
since it amounts to replacing $G$ by an automorphism group of a variety.
One of the main points of  
\cite{twcv} was to understand the definition of 
$\cI$-graded $G$-local systems,  as an action of a certain local system $\cT$  of  infinite dimensional tori (its fibres are isomorphic to 
the exponential torus of \cite{MR91}).
Thus in the nonlinear case a Stokes graded local system
will involve a nonlinear local system $V$ with a 
preferred locally constant torus action in sectors at 
infinity.
Generically one would expect this torus to have dense orbits, 
so the fibres of $V$ on sectors will have dense tori in them.
For the local systems formed by the wild character varieties (cf. \cite{smid, gbs}), 
one expects to be able to write down these tori explicitly, 
analogously to the well known explicit description of the (tame) nonlinear monodromy (this question was raised in 
\cite{ramispaul}).
In some examples such tori are evident by combining the complex WKB viewpoint \cite{voros}
and the TQFT approach to 
meromorphic connections (initiated in \cite{smid, saqh}
and completed in \cite{fission, gbs, twcv}):
Each generic Stokes graph (\cite{voros} p.271)
corresponds to dividing the surface into pieces, and approximating the connection by the Airy equation in each piece. 
Translating to the TQFT approach, this amounts to fusing together several copies of the fission space 
$\cB$ corresponding to the Airy equation. But the fission space of the Airy equation is just a copy of $\IC^*$, and fusion means taking the product of such spaces, yielding a torus $(\IC^*)^n$
(see the end of \cite{ihestalk2017} for the pictures explaining this).

\subsection{Layout of the article}

To orient the reader,
\S\ref{sn: onepagesummary}
gives a short (one page) summary of 
the data canonically attached to a connection
that is going to be studied here.
The subsequent sections then recall basic notions, 
related to linear algebra \S\ref{sn: linalg}
(gradings, filtrations, splittings, relative position) and to local systems/covers 
\S\ref{sn: topo}.
Next \S\ref{sn: irclass} recalls the notion of an 
irregular class (in the general linear case) and the
resulting 
notion of irregular curve/wild Riemann surface
(i.e. a curve with some marked points, each 
equipped with an irregular class).
The core of the article consists of 
\S\S\ref{sn: SFLS},\ref{sn: SGLS},\ref{sn: SLS}
that define
Stokes filtered local systems, Stokes graded local systems 
and Stokes local systems, respectively.
Some basic properties are then established leading up to 
the proof of  Thm. \ref{thm: mainresult} in 
\S\ref{sn: canl splittings}.
Next \S\ref{sn: wcvs} reviews the implications of this for the wild character varieties, and
explains the resulting notion of 
wild nonabelian periods/wild Wilson loops.
For completeness Apx. \ref{sn: anbbx} summarises the analytic results (presented as black boxes) 
needed to 
attach such topological data to a connection. 
In the final section some of the basic ideas used by 
Stokes to get to this picture are sketched. 

Note that the word ``Stokes'' will often be used
to indicate possible discontinuity: 
whereas a graded local system is continuous, 
the graded pieces of a ``Stokes graded local system'' may have discontinuities at the singular directions, 
and similarly 
the filtered pieces of a 
``Stokes filtered local system'' may have discontinuities at the oscillating directions.
This is in the spirit of \cite{stokes1857}.
It is worth emphasising that there are thus {\em two} types of discontinuity that occur, in general in different directions. 
(On the other hand a Stokes local system has no discontinuities, although it lives on the surface obtained by removing the tangential punctures
and is only graded in the halos.)

\begin{figure}[h]
\centerline{
  \resizebox{5cm}{!}{
    \begin{tikzpicture}[scale=0.5, level/.style={thick}]
  \draw[thick, rotate=0] (0,0) ellipse (5cm and 2cm);
  \draw[thick, rotate=60] (0,0) ellipse (5cm and 2cm);
  \draw[thick, rotate=120] (0,0) ellipse (5cm and 2cm);
    \end{tikzpicture}
  }
} 
\caption{Example Stokes diagram, see \S\ref{ssn: weber}. }%
\label{fig: atomic}
\end{figure}
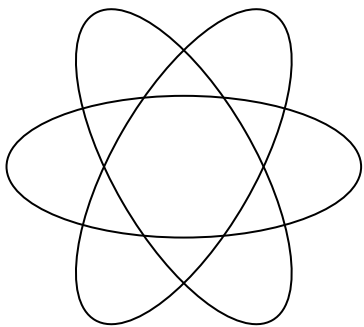

\newpage

\section{Summary of some data canonically determined
by a connection}\label{sn: onepagesummary}

Let $\Si$ be a smooth compact complex algebraic curve, $\ba\subset \Si$ a finite subset, and 
$\Si^\circ=\Si\setminus \ba$.
Let $\pi:\wh \Si\to \Si$ be the real oriented blow up 
of $\Si$ at $\ba$ and let 
$\partial=\pi^{-1}(\ba)$.
Thus $\partial=\bigcup_{a\in \ba}\partial_a$ is a 
collection of circles, and a point $d\in \partial_a$ 
is a real oriented direction in the tangent space
$T_a\Si$ at $a\in \ba\subset \Si$.
There is a canonically defined covering space $\cI\to \partial$, the 
{\em exponential local system} (see \S\ref{ssn: exp loc sys} below).

Suppose $(E,\nabla)$ is a connection on an algebraic vector bundle on $\Si^\circ$.
Then the following data are canonically defined:

1) The local system $V=\ker(\nabla^{an})$ of analytic solutions. It is a locally constant sheaf of finite dimensional complex vector spaces on $\Si^\circ$, and 
we extend it to $\wh \Si$ in the obvious way.

2) An irregular class $\Th:\cI\to \IN$. This will be recalled in detail below---in brief it amounts to choosing the exponential factors plus their multiplicities at each marked point.
The  exponential factors will be viewed as a finite covering space $I\to \partial$.

3) The singular directions $\IA\subset \partial$
and the Stokes (oscillating) directions
$\IS\subset \partial$. They are both finite sets.

4) A total ordering of the finite set $I_d$ (the fibre of $I\to \partial$ at $d$) for any non-oscillating direction 
$d\in \partial\setminus \IS$. This is the natural
dominance ordering of the exponential factors.

5) For any $d\in \partial\setminus \IS$, a 
filtration of $V_d$ by the finite set $I_d$,
the {\em Stokes filtrations}.

6) For any $d\in \partial\setminus \IA$, a 
grading of $V_d$ by the finite set $I_d$,
the {\em Stokes decompositions} or 
{\em Stokes gradings}. 

7) An $I$-graded local system $V^0\to \partial$ of the same rank as $V$, the {\em formal local system} (beware it is not a grading of $V$).

8) For any $d\in \partial\setminus \IA$, a 
linear isomorphism $\Phi:V_d^0\to V_d$, the 
{\em gluing maps}.

9) For any $d\in\IA$, two unipotent 
automorphisms: $g_d\in \GL(V_d)$,
the wild monodromy automorphism, and
$\IS_d\in \GL(V^0_d)$ the Stokes automorphism.

10) A local system $\IV\to \wt \Si$, the {\em Stokes local system}, where $\wt \Si\subset \wh \Si$
is the auxiliary surface obtained by removing a tangential puncture 
$e(d)\in \Si^\circ$ near each $d\in \IA$.

These data  have lots of properties not mentioned here and are not independent (in particular the wild monodromy and the Stokes automorphism are essentially the same thing). 
Various subsets of these data can be (and have been) precisely axiomatised to encode the category of connections.
Our basic aim is to describe this story and the relations between these data.
Note this list is not comprehensive (although any other data will be a function of the data here)---the sequel will discuss the {\em Birkhoff gradings}.

\newpage

\section{Linear algebra}\label{sn: linalg}

\subsection{Gradings}

Let $V$ be a finite dimensional complex vector space and let $I$ be a set.
A {\em grading} $\Ga$ of $V$ by $I$ is 
the choice of a subspace $\Ga_i=\Ga(i)\subset V$ 
for each $i\in I$
so that there is
a direct sum decomposition 
$$V=\bigoplus_{i\in I} \Ga(i).$$
Some of the $\Ga(i)$ are allowed to be the zero subspace, so $I$ could be much bigger than the dimension of $V$.
An {\em $I$-graded vector space} is a pair $(V,\Ga)$.

The multiplicity (or {\em dimension vector})
of a grading 
is the element $\Th\in \IN^I$ with components
$\Th(i) = \dim(\Ga_i)$.
In other words it is the map 
$\Th:I\to \IN; i\mapsto \dim(\Ga_i)$.
Given $(V,\Ga)$ an index $i\in I$ is {\em active} if
$\Ga_i$ is nonzero, i.e. $\Th(i)\neq 0$.

Let $\Aut(V,\Ga)\subset \Aut(V)=\GL(V)$ denote the group of graded automorphisms, i.e. $g\in \GL(V)$ such that 
$g(\Ga_i)=\Ga_i$ for all $i\in I$.

A grading is {\em full} or {\em toral} if $\dim(\Ga_i)\le 1$ for all $i$, so that $\Aut(V,\Ga)$ is a torus.
A basis $\{ e_i\}$ of $V$ determines a full grading by taking 
$\Ga_i$ to be the line $\IC e_i\subset V$.
Conversely a full grading determines a basis up to the action 
of a torus (the choice of a basis of each one dimensional subspace $\Ga_i$).

Sometimes (if $V$ just has one grading) 
the graded pieces will be written
$V(i)=V_i=\Ga_i$. 
Then $V$ is said to be an $I$-graded vector space, and  $\GrAut(V)\subset \GL(V)$ will denote the group of 
graded automorphisms.

\subsection{Filtrations}

Now suppose the set $I$ is given a total ordering $\le$.
An {\em $I$-filtered vector space} is a pair $(V,F)$
where $V$ is a complex vector space and $F$ is a filtration of $V$ indexed by $I$, i.e. 
a collection of subspaces 
$F_i\subset V$ for each $i\in I$ such that if $i\le j$ then $F_i \subset F_j$.
We will sometimes write $F(<\!i)=\sum_{j<i}F_j$, and $F(i)=F_i$.
Here $j<i$ means that $j\le i$ and $j\neq i$
(with this understood the order is determined equivalently by the 
binary relation $<$ or $\le$).
A map between $I$-filtered vector spaces $(V,F),(W,G)$
is a linear map $\varphi:V\to W$ 
such that $\varphi(F_i)\subset G_i$ for all $i\in I$.

The {\em associated graded vector space} of $(V,F)$  is 
the $I$-graded vector space with graded pieces
$\Gr_i(V,F) := F(i)/F(<\!i)$, i.e. it is the external direct sum
$$\Gr(V,F):=\bigoplus_I \Gr_i(V,F),\qquad \Gr_i(V,F) := F(i)/F(<\!i).$$ %
A map $\varphi:(V,F)\to(W,G)$ 
of $I$-filtered vector spaces 
induces a map $\Gr(\varphi):\Gr(V,F)\to\Gr(W,G)$
of graded vector spaces. 
By definition the multiplicity (or {\em dimension vector})
$\Th:I\to \IN$
of a filtered vector space is that of its associated graded.
Thus given $(V,F)$ then $i\in I$ will be said to be active if
$\Gr_i(F)$ is nonzero.

Recall that a flag (of subspaces of $V$) is a nested collection of distinct
subspaces $0=V_0\subset V_1\subset V_2\subset \cdots V_{k-1}\subset V_k =V$
for some $k$.
It is {\em full} if $\dim(V_i)=i$.
Thus an $I$-filtration consists of the choice of 
 a flag plus the choice of 
a labelling of the subspaces by elements of $I$.
If the dimension vector (and the ordering of $I$) is fixed 
then there is unique labelling, so choosing a filtration with a given dimension vector is the same as choosing a flag (with the subspaces of the right dimensions).

Since $I$ is ordered, 
an $I$-grading $\Ga$ of $V$ determines a filtration 
$\cF=\cF(\Ga,\le)$ of $V$ defined by:
$$\cF_i := \bigoplus_{k\le i} \Ga(k).$$
This is the  ``associated filtration'' determined by the ordering.

\subsection{Splittings}

A {\em splitting} of an $I$-filtered vector space $(V,F)$ 
is the choice of an $I$-grading $\Ga$ 
of $V$
such that $F$ equals the associated filtration 
$\cF(\Ga,\le)$.

It will be useful to think of splittings in terms of 
isomorphisms with the associated graded, as follows.
Let $(V,F)$ be an $I$-filtered vector space.
Let $\IF$ denote the associated filtration of $\Gr(V,F)$
so that  $\IF_i = \bigoplus_{j\le i} \Gr_j(V,F)$. 
Clearly the filtration $\IF$ has a preferred splitting, so we can canonically identify its associated graded with $\Gr(V,F)$.
\begin{lem}\label{lem: splitting}
Giving a splitting of $(V,F)$ is the same as giving 
an  isomorphism $$\Phi:(\Gr(V,F),\IF)\mapright{\cong} (V,F)$$ 
of filtered vector spaces
such that 
the associated graded map 
$$\Gr(\Phi):\Gr(V,F)\to\Gr(V,F)$$
is the identity, i.e. the map 
$\Gr_i(F)\mapright{\Phi} F_i\onto \Gr_i(F)$
is the identity for each $i\in I$.
\end{lem}
\pf
Straightforward. Given such $\Phi$, the splitting is given by 
$\Ga_i=\Phi(\Gr_i(V,F))$.
Conversely, given a splitting, the map $\Ga_i\into F_i\onto \Gr_i(V,F)$ is an isomorphism and its inverse gives the $i$-component of $\Phi$.
\epf

This will be crucial to define the wild monodromy
automorphism
determined by a pair of compatible gradings, and motivate the notion of Stokes local system.

\subsection{Wild monodromy---Relative positions of pairs of gradings}

Given two bases of a vector space $V$ (indexed by the same set) 
it is obvious that
there is a unique automorphism of $V$ taking one basis to the other.
However if only the underlying gradings are given 
(and not the bases
themselves) then in general there is not a preferred 
automorphism  taking one grading to the other.
However it turns out that,
for a  special class of pairs of gradings, 
there is indeed a preferred 
automorphism taking one grading to the other (uniquely determined by the gradings).

Suppose $I$ is a set and $V$ is a vector space.
Two $I$-gradings $\Ga_1,\Ga_2$ of $V$ are {\em compatible}
if there is some ordering $\le$ of $I$
such that the associated 
filtrations are equal:
\beq\label{eq: com filt}
\cF(\Ga_1,\le) = \cF(\Ga_2,\le).
\eeq
Thus this common filtration
is split by both gradings.
Note this implies both gradings have the same dimension vector.
Also note that if $\Ga_1,\Ga_2$ are compatible, then there may well be several different orderings for which 
their filtrations are equal.
The aim  here is to show there is a preferred automorphism relating any two compatible gradings.

Indeed if $F$ denotes the common filtration
\eqref{eq: com filt}
then Lemma \ref{lem: splitting} implies there are distinguished linear isomorphisms
$$\Phi_1,\Phi_2: \Gr(V,F)\to V$$
determined by the splittings $\Ga_1,\Ga_2$ respectively, 
so there is a preferred automorphism
$$g=g(\Ga_1,\Ga_2) = \Phi_2\circ \Phi_1^{-1}\in \GL(V)$$
taking $\Ga_1$ to $\Ga_2$, the {\em wild monodromy}.
Of course since the gradings are given, it is 
simpler to work with the inverse 
$\Psi_i:=\Phi_i^{-1}$,
which is  just the natural graded isomorphism 
$$\Psi_i:(V,\Ga_i)\ \mapright{\cong}\  \Gr(V,F),\qquad
\text{so that } g= \Psi_2^{-1}\circ \Psi_1.$$

An analysis of the set of splittings of a given filtration enables to see it is independent of
the choice of $F$.
Fix an ordering $\le$ of $I$ and an
$I$-filtration $F$ of $V$.
Let $\Th$ be its dimension vector.
Let $\cG$ be the set of all $I$-gradings of $V$
of dimension $\Th$, 
and let 
$\Splits(F)\subset \cG$ 
be the set of all $I$-gradings 
that split $F$, i.e. 
the gradings $\Ga$ such that $\cF(\Ga,\le)=F$.

The group $G=\GL(V)$ acts transitively on $\cG$, so $\cG$  is a homogeneous space.
Given a basepoint $\Ga_1\in \cG$
then the map 
$G\onto \cG; g\mapsto g(\Ga_1)$ identifies $\cG\cong G/H$
where $H=\GrAut(V,\Ga_1)\subset G$.

In this way the subset $\Splits(F)\subset \cG$ 
corresponds to the subset $P/H\subset G/H$ where 
$P=P(F)=\Aut(V,F)\subset G$ is the group of filtered automorphisms
(i.e. the parabolic subgroup fixing the flag underlying $F$).

Now the quotient $P\to P/H$ has a natural slice: 
the unipotent radical $U=U(F)=\Rad_u(P)\subset P$ 
is the maximal unipotent normal subgroup of $P$, 
and it has the property that the product map
$U\times H \to P$ is an isomorphism of varieties 
($P$ is isomorphic as a group to the semidirect product $H \sdp U$)\footnote{In concrete terms, choosing suitable bases, $P$ is a block upper triangular subgroup of $G$, 
$H$ is the block diagonal subgroup of $P$,
and $U$ is the subgroup of $P$ with an identity matrix in each diagonal block.}.
Thus the above map $G\onto \cG$ restricts to an isomorphism 
$U(F)\cong \Splits(F)$, mapping $U$ isomorphically 
onto $\Splits(F)$.

More intrinsically this shows that $\Splits(F)$ is a torsor
(principal homogeneous space) for $U(F)$---the 
natural action of $U(F)$ on $\Splits(F)$ 
is free and transitive.
In particular for any two elements 
$\Ga_1,\Ga_2\in\Splits(F)$ 
there is a unique element 
$g=g(\Ga_1,\Ga_2)\in U(F)$ such that 
$g(\Ga_1)=\Ga_2$ (it equals the wild monodromy
as that has this property).
This element $g\in \GL(V)$ is in fact 
uniquely determined by the pair $\Ga_1,\Ga_2$
and does not depend on the choice of filtration $F$ that 
they both split:

\begin{lem}
If $\Ga_1,\Ga_2$ are two compatible gradings then
there is a  preferred 
unipotent element
$g=g(\Ga_1,\Ga_2)\in \GL(V)$
taking $\Ga_1$ to $\Ga_2$, 
the wild monodromy automorphism.
It is the unique element $g\in\GL(V)$ such that 
$g(\Ga_1)=\Ga_2$ and 
$\Gr(g)=1\in \GrAut(\Gr(V,F))$, for any filtration $F$ split by both gradings.
\end{lem}

\pf
It just remains to show that  if the gradings give the
same filtration for two different orderings then
the wild monodromy elements are the same: 
Let $F_1,F_2$ be the two filtrations (for two different orderings), and suppose 
$$\Ga_1,\Ga_2\in \Splits_{12}:=\Splits(F_1)\cap \Splits(F_2).$$
The result then follows from the fact that 
$\Splits_{12}$ is a torsor for the 
group $U(F_1)\cap U(F_2)$,
which is the unipotent radical of 
$P(F_1)\cap P(F_2)$.
\epf

Note in general (in the non-toral case) 
the wild monodromy is not the only unipotent element
taking one grading to the other.
In practice extra choices are often made to give bases---the discussion here shows in general that the resulting 
automorphism only depends on the gradings.
The following easy lemma will be useful.

\begin{lem}\label{lem: monod invt vectors}
Suppose $\Ga_1,\Ga_2$ are compatible $I$-gradings of $V$
with wild monodromy $g=g(\Ga_1,\Ga_2)\in \GL(V)$. 
If $v\in \Ga_1(i)\cap\Ga_2(i)$ for some $i\in I$
then $g(v)=v$.
\end{lem}
\pf It is clear that $\Psi_1(v)=\Psi_2(v)$, so the result follows. \epf

\begin{eg}
If $\dim(V)=\#I=2$ and $\cG$ is the set of 
full $I$-gradings of $V$, then 
$\cG$ is just 
the space of injective maps $\Ga:I\into \IP(V)$, since any two distinct lines in $V$ are linearly independent.
In other words it is the space of pairs of distinct points of the sphere $\IP(V)$, labelled by $I$.
If $I=\{ 1,x\}$  then 
two gradings $\Ga_1,\Ga_2$ are compatible if and only if 
either 
$\Ga_1(1)= \Ga_2(1)$ or  $\Ga_1(x)= \Ga_2(x)$.
If $F$ is the filtration determined by $\Ga_1$ with the ordering $1<x$, then the filtration is just $F(1)=\Ga_1(1)\subset F(x)=V$.
Thus $\Splits(F)$ is the space of maps 
$\Ga:I \into \IP(V)$ such that 
$\Ga(1)=F(1)$. This amounts to the choice of a point 
$$ \Ga(x)\in \IP(V)\setminus \{F(1)\}$$
of the affine line given by the sphere punctured at $F(1)$.
Thus $\Splits(F)= 
\IP(V)\setminus \{F(1)\}\cong \IA^1$.
It is a torsor for the group $U(F)\cong \{\bsmx 1 & * \\ 0 & 1 \esmx\}\cong \IA^1$. 
\end{eg}

\subsection{Median gradings} \label{sn: median gr}
If $\Ga_0,\Ga_1$ are two compatible
$I$-gradings of $V$ 
let $g=g(\Ga_0,\Ga_1)\in \GL(V)$
be the wild monodromy
relating them, so that 
$\Ga_{1} = {g}(\Ga_0)$.
Since $g$ is unipotent it has a unique unipotent square root $\sqrt{g}\in \GL(V)$
and so the {\em median grading}
$$\Ga_{1/2} := \sqrt{g}(\Ga_0)$$
is well defined.
Of course $\sqrt{g} = g(\Ga_0,\Ga_{1/2}) = 
g(\Ga_{1/2},\Ga_1)$.
More generally 
there is a canonically 
defined path of gradings $\{\Ga_t\st t\in [0,1]\}$
connecting $\Ga_0$ to $\Ga_1$ given by 
$\Ga_t := g_t(\Ga_0)$ where 
$g_t=\exp(tX)$ for the 
unique nilpotent logarithm $X$ of $g$.

In the example above with $\dim(V)=2$,
the median grading between two points of $\cG(F)\cong \IA^1$ 
is just the midpoint of the real line segment in $\IA^1$ between them.  

This appears classically  as the way to get real bases/gradings of a real differential equation 
(cf. \cite{dinglebook} p.8, \cite{MR91} p.358).
The main statement boils down to  the following.
Suppose $I$ is a set, 
$V_\IR$ is a real vector space and $V=\IC\otimes V_\IR$ is its complexification.
\begin{lem}
If $\Ga_0,\Ga_1$ are compatible $I$-gradings of $V$
such that $\Ga_0=\bar{\Ga_1}$, then
$\Ga_{1/2}$ is {\em real}.
\end{lem}
\pf
Write $g=g(\Ga_0,\Ga_1), s=\sqrt{g}$. 
Note 
$\Ga_0=\bar{g}\bar{\Ga_0}=\bar{g}\Ga_1$
and $\Ga_0=g^{-1}\Ga_1$, so by uniqueness
$\bar{g}=g^{-1}$, so that $g=\exp(X)$ for a unique 
nilpotent $X$ such that $\bar{X}=-X$.
Thus $s=\exp(X/2)$ satisfies $\bar{s}=s^{-1}$ and the result follows:
$$\bar{\Ga_{1/2}} = \bar{s}\bar{\Ga_0}=s^{-1}\Ga_1
=s^{-1}g\Ga_0= \Ga_{1/2}.$$
\epf

\subsection{Bounding the  wild monodromy}

Suppose $I$ is a set 
equipped with a partial order $\prec$.
Thus $\prec$ is a subset of  $I\times I$
satisfying various axioms.
The pair $(I,\prec)$ can be viewed as a quiver
where a relation $i\prec j$ corresponds 
to an arrow $i\from j$.
The choice of $\prec$ can be used to restrict the 
wild monodromy, as follows.

Recall that a total order $<$ on $I$ is said to ``extend the partial order 
$\prec$'', if $\prec$ is a subset of $<$ (they are both subsets of $I\times I$).

\begin{defn}
The wild monodromy 
of a pair of compatible $I$-gradings 
$\Ga_1,\Ga_2$
of $V$ is  ``bounded by $\prec$'', or  ``satisfies the Stokes conditions'', if the associated filtrations
\beq\label{eq: Gstokescondns}
\cF(\Ga_1,<)=\cF(\Ga_2,<) 
\eeq
are equal, for any total order $<$ extending $\prec$.
\end{defn}

\begin{rmk}\label{rmk: empty prec}
For example if $\prec$ is empty this just means that $\Ga_1=\Ga_2$.
\end{rmk}

The Stokes condition can be reformulated in terms of Stokes groups as follows.
For $k=1,2$ let
 $\ISto_k\subset \GL(V)$ be the connected 
 unipotent group 
with Lie algebra
\beq\label{eq: st gp}
\Lie(\ISto_k) = \bigoplus_{i\prec j} 
 \Hom(\Ga_k(j),\Ga_k(i))\subset \End(V).
\eeq
Note $\Lie(\ISto_k)$ is the space of representations of the quiver $(I,\prec)$ on $(V,\Ga_k)$\footnote{The 
general fact used here is that 
if $V$ is 
$I$-graded then a partial order $\prec$ on 
$I$ determines a unipotent subgroup $\ISto_\prec$ 
of $\GL(V)$ whose Lie algebra is the space of representations of the quiver 
$(I,\prec)$ on $V$.}.
Consider the set 
$\GrIso_{12}=\{ g\in \GL(V)\st g(\Ga_1(i))=\Ga_2(i) \text{ for all $i\in I$} \}$
of graded isomorphisms
from $(V,\Ga_1,)$ to $(V,\Ga_2)$.

\begin{lem}\label{lem: st.condns}
The following conditions are equivalent:

1) $\GrIso_{12}\cap \ISto_1\neq  \varnothing$, \ \ 
2) $\GrIso_{12}\cap \ISto_2\neq  \varnothing$,

3) $g(\Ga_1,\Ga_2)\in \ISto_1,$ \ \ \ \ \ \ 
4) $g(\Ga_1,\Ga_2)\in \ISto_2,$ 

5) The gradings $\Ga_1,\Ga_2$ satisfy the Stokes conditions determined by $\prec$.

\noindent
If so, $\ISto_1=\ISto_2$ 
and both sets 1),2) 
 contain exactly one point $g(\Ga_1,\Ga_2)\in \GL(V)$.
\end{lem}

\pf
If $\Ga_1$ is fixed then any other
$I$-grading $\Ga_2$ of $V$ 
(with the same dimensions as $\Ga_1$)
can be specified by choosing an element
$g\in\GL(V)$ and defining
$\Ga_2(i)=g(\Ga_1(i))$.
Then $g\in \GrIso_{12}$.
Moreover the two Stokes groups in $\GL(V)$ determined by the two gradings are then conjugate by $g$: 
$\ISto_2= g\circ \ISto_1 \circ g^{-1}.$ 
Now suppose we fix both gradings and assume $g\in  \GrIso_{12}$. 
Then it is immediate that 
$g\in  \ISto_1$ if and only if 
$g\in\ISto_2=g\ISto_1 g^{-1}$, proving 1) and 2) are equivalent.
Moreover if $g\in \ISto_1$ then clearly
$\ISto_2=g\ISto_1g^{-1}=\ISto_1.$
If moreover $g_1\in \GrIso_{12}\cap \ISto_1$
then $g^{-1}g_1\in \ISto_1\cap \Aut(V,\Ga_1)=\{1\}$ so $g_1=g$.
This shows the last statement holds and 
establishes the equivalence with 3,4).
Finally to see 1-4) are equivalent to 5)  observe that  the Stokes group is the unipotent radical of the intersection of 
the parabolic subgroups preserving 
$\cF(\Ga_1,<)$, as $<$ varies.
\epf

\subsection{Stokes conditions on pairs of filtrations}

The above discussion will be used to control the 
jumps of the Stokes gradings across $\IA$.
The jumps of the Stokes filtrations 
across $\IS$ are controlled as follows.  

Let $I$ be a set equipped with two orders $<_1,<_2$.
Write $I_1=(I,<_1), I_2=(I,<_2)$ for the corresponding
ordered sets.
Fix a complex vector space $V$ and let
 $F_1$ be an $I_1$-filtration of $V$, and let
 $F_2$ be an $I_2$-filtration of $V$, both with the same 
 dimension vector $\Th:I\to \IN$.

\begin{defn}
The pair of filtrations $(F_1,F_2)$ 
``satisfy the Stokes conditions'' if there 
is an $I$-grading $\Ga$ of $V$ of dimension $\Th$ such that
\beq\label{eq: Fstokescondns}
F_1=\cF(\Ga,<_1)\quad \text{and}\quad  F_2=\cF(\Ga,<_2).
\eeq
\end{defn}

The notion of relative position of pairs of flags has been much studied (cf. \cite{del-lusz76} p.116) and one can show that the Stokes conditions above are the same as fixing the relative position of the pair of flags
to be that determined by the pair of orderings 
(although, beyond the toral case, not every possible relative position will arise in the Stokes setting).

\section{Topological basics}\label{sn: topo}

\subsection{Local systems, transport, monodromy}

Recall that if $M$ is a connected topological manifold 
then a local system of sets 
$I$ on $M$ is a locally constant sheaf of sets, 
and that this is the same thing as (the sections of) a covering space $I\to M$.
Given an open cover of $M$ then $I$ 
can be described in terms of
{\em constant} clutching maps.

A local system of vector spaces is a local system $V\to M$ 
for which the fibres are 
vector spaces and this structure is preserved (the clutching maps are constant linear isomorphisms).

Given  a local system $V\to M$ then 
a path $\ga$ in $M$ from $p$ to $q$ 
determines an isomorphism 
$$\rho_\ga=\rho(\ga):V_p\to V_q$$
between the corresponding fibres, the parallel transport map.
(It is defined since the path has a unique lift to the covering
once the initial point in $V_p$ is specified.)
Homotopic paths give the same map.
In particular if a  basepoint $b\in M$ is fixed 
then the monodromy representation 
$$\rho:\pi_1(M,b)\to \Aut(V_b)$$
is defined by transporting 
 points of the fibre $V_b$ around loops. 
Two local systems $V,V'$ 
are isomorphic if and only if there is an isomorphism
$V_b\cong V'_b$ intertwining their monodromy representations.

If $V$ is a local system of rank $n$ vector spaces then a framing of $V$ at $b\in M$ 
is a basis of the fibre of $V$ at $b$, i.e. 
an isomorphism $\phi:\IC^n\to V_b$.
Given a framing,
the monodromy representation 
can be viewed as taking values in $\GL_n(\IC)$.

\subsection{Graded local systems}\ 
If $I\to M$ is a fixed covering space, then an 
``$I$-graded local system'' (of vector spaces) 
is  a local system $V\to M$ of vector spaces together with 
a pointwise grading
$$V_p = \bigoplus_{i\in I_p} V_p(i)$$
of the fibres of $V$ by the fibres of $I$,
such that the grading is locally constant, i.e.
$$\rho_\ga(V_p(i)) = V_q(\rho_\ga(i))$$
for any path $\ga$ from $p$ to $q$
(where $\rho_\ga(i)\in I_q$ is the parallel transport of $i\in I_p$).
This is the same thing as a local system on the covering space $I$, but viewed from $M$.

\subsection{Real oriented blow-ups and tangential basepoints}\label{ssn: robups}

If $\Si$ is a smooth complex algebraic curve
and $a\in \Si$,
the real oriented blow-up of $\Si$ at $a$
is the surface with boundary 
$$\pi:\wh \Si\to \Si$$
obtained by replacing the point $a$ by the circle
$$\partial = (T_a\Si\setminus\{0\})/\IR_{>0}$$
of real oriented directions emanating from 
$a$. Here $T_a\Si$ is the tangent space.
Thus $\pi$ restricts to an isomorphism 
$\wh\Si\setminus \partial\cong \Si^\circ := 
\Si\setminus \{a\}$ onto the punctured surface.

A useful construction is the following:
If $V\to \Si^\circ$ is a local system of $n$-dimensional vector spaces
and $d\in \partial$, then 
there is a well defined complex vector space
$V_d$ of dimension $n$,
the 
``fibre of $V$ at the tangential basepoint $d$''.
One approach is to note that the inclusion
$\Si^\circ \to \wh \Si$ is a homotopy equivalence
so restricting local systems on $\wh \Si$ to 
$\Si^\circ$ gives an equivalence of categories.
(Any covering of $\Si^\circ$ extends 
uniquely to $\wh \Si$.)
In this way local systems on $\Si^\circ$ and on $\wh \Si$ will 
henceforth be viewed as the same thing.
Then $V_d$ is the fibre of $V\to \wh \Si$ at the tangential basepoint $d\in \partial\subset \wh \Si$. 
More concretely
$V_d$  is 
the space of sections of $V$ on 
germs of open sectors (in $\Si^\circ$)
at $a$ containing the 
direction $d$ (cf. \cite{deligne-3pts} p.85, 
\cite{Mal-irregENSlong} p.386).

\subsection{Extended intervals/sectors}\label{ext ints}

Suppose $a\in \Si$ and $\pi:\wh \Si\to \Si$ is the real 
oriented blow-up at $a$, and $\partial$ is the circle 
$\pi^{-1}(a)\subset \wh \Si$.
This section will set up notation for 
universal covers $\wt \partial$ of $\partial$ and
intrinsically 
extending a local system on $\partial$ to a sector/interval in $\wt \partial$.

1) The notion of angle is well defined so there is an intrinsic action of $\IR$ on $\partial$. 
Thus if $d\in \partial$ and 
$\al\in \IR$ then $d\pm\al$ are  well-defined points of 
$\partial$.

2) The choice of a point $d\in \partial$ determines a universal covering $\wt \partial$ of $\partial$, namely the homotopy classes of paths starting at $d$.
It comes equipped with a covering map 
$\wt \partial\onto \partial$ 
(taking a path to its endpoint) 
and with 
a point $\wt d\in \wt \partial$ lying over $d$ (i.e. the trivial path).

Now, for positive $\al,\be\in \IR$ define 
$\Sect_d(-\al,\be)$ 
to be the interval 
$(\wt d-\al,\wt d+\be)\subset \wt \partial$ in the universal cover determined by $d$.
Similarly $\Sect_d[-\al,\be) =
[\wt d-\al,\wt d+\be)$ etc. 
If $U\subset \partial$ is an interval, define 
\beq\label{eq: sectU}
\Sect_U(-\al,\be) = \bigcup_{d\in U}\Sect_d(-\al,\be).
\eeq

3) A local system $V\to \partial$ canonically determines a local system
$\wt V\to \Sect_d(-\al,\be)$, 
for any  positive 
$\al,\be\in \IR$, and $V_d=\wt V_{\wt d}$. 
It is defined as the pullback of $V$ along the (\'etale) map 
$\Sect_d(-\al,\be)\to \partial$.
For small $\al,\be$ this is just the restriction to 
$(d-\al,d+\be)\subset \partial$.
There is a canonical bijection between $V_d$ and the 
space of sections of $\wt V$ on $\Sect_d(-\al,\be)$
for any $\al,\be$ 
(taking a section to its value in $\wt V_{\wt d} = V_d$).

\section{Irregular classes and associated topological data}\label{sn: irclass}
Fix $\ba\subset \Si$ 
as in \S\ref{sn: onepagesummary}.
Let $\wh \Si\to \Si$ be the real oriented blow-up,
with boundary $\partial$ (a finite set of circles), cf. \S\ref{ssn: robups}.
The aim of this section is to recall the notion of irregular class 
(from \cite{twcv}) and to show (following mainly 
\cite{deligne78in, Mal-irregENSlong, BJLproper, MR91, L-R94}) that this determines the following 
(topological) data:

1) an integer $n$ (the rank),

2) a finite cover $I\to \partial$, 

3) a dimension vector for $I$, i.e. a map $\Th:I\to \IN_{>0}$ 
constant on each component circle, 
such that
$\sum_{i\in I_d}\Th(i)=n$ for all $d\in \partial$, 

3) two finite sets $\IA,\IS\subset \partial$

4) A total order $<_d$ of $I_d$ for each 
$d\in \partial\setminus \IS$ (constant as $d$ moves in each component of this set)

5) A partial order  $\prec_d$ of $I_d$ for each 
singular direction $d\in \IA$, the Stokes arrows.

This data alone will be sufficient to define the topological data classifying connections 
in the next three sections, 
so the topologically minded reader could (in the first instance) skip the following discussion
and jump to \S\ref{sn: SFLS}. 

\subsection{Exponential local system}
\label{ssn: exp loc sys}

The exponential local system is a natural 
covering space (local system of sets) 
$\pi:\cI\to \partial$.
To simplify notation suppose $\ba=\{a\}$ 
is just one point so 
$\partial$ is a single circle (the extension to multiple points is immediate).
If $z$ is a local coordinate on $\Si$ 
vanishing at $a$ then
local sections of $\cI$ over open subsets of 
$\partial$ are functions that may be written as finite sums of the form 
\beq\label{eq: q as sum}
q= \sum a_i z^{-k_i}
\eeq
where  $a_i\in \IC$, and $k_i\in \IQ_{>0}$.
An intrinsic (coordinate independent) 
construction of $\cI$ is given in \cite{twcv} Rmk 3
(in which case sections of $\cI$ are certain equivalence classes of functions, but that will make no difference in the use of these functions below).
Thus $\cI$ is the disjoint union of a vast collection of circles $\<q\>$, each of which is a finite cover of $\partial$. A local function $q$  determines the circle $\<q\>$ by analytic continuation around $\partial$.
Thus algebraically $\<q\>$ 
encodes the Galois orbit of $q$.
Let $\ram(q)$ denote the degree of the cover 
$\pi:\<q\>\to \partial$ (the lowest common multiple of the denominators of the $k_i$ present in the expression for $q$).
The {\em slope} of $\<q\>$ is the largest $k_i$
occuring in \eqref{eq: q as sum}.
The {\em tame circle} is the circle $\<0\>\subset \cI$.
The functions $q$ occur as the exponents of the exponential factors $e^q$ that occur in local solutions of meromorphic differential equations (hence the name ``exponential local system''). An isomorphic local system $d\cI$ (whose sections are one-forms) was used in \cite{deligne78in, Mal-irregENSlong}.

\subsection{Irregular classes}

An irregular class is a map 
(a dimension vector) 
$\Th:\cI\to \IN$, 
assigning an integer to each component of $\cI$, equal to zero for all but a  finite number of circles.
(Thus $\Th$ is constant on each component circle,
so amounts to a map $\pi_0(\cI)\to \IN$.)
The {\em rank} of an irregular class is the integer
\beq 
n = \rk(\Th)= \sum_{i\in \cI_d}\Th(i)\in \IN
\eeq
for any $d\in \partial$.
(If $\partial$ has several components then $\Th$ should be such that the rank is the same for $d$ in any component.)

A {\em finite subcover} is a subset $I\subset \cI$
such that $\pi:I\to \partial$ is a finite cover.
An irregular class determines a finite subcover, 
the {\em active exponents} $I=\Th^{-1}(\IN_{>0})$.
Thus an irregular class is a 
finite subcover plus a positive integer for each component (we will often omit to write the 
integers and say that $I$ is an irregular class).

Any $\cI$-graded local system $V\to \partial$ of vector spaces has an irregular class (taking the dimensions of the graded pieces).
It will thus follow that any connection determines an irregular class 
(the associated graded of its Stokes filtration
is an $\cI$-graded local system).

Given irregular classes $I_1,I_2$ then there are well-defined irregular classes 
$I_1^\vee$, $\End(I_1)$, $I_1\otimes I_2$, 
$\Hom(I_1,I_2)=I_2\otimes I_1^\vee$.
In brief if $V_k\to \partial$ is any $I_k$-graded
local system for $k=1,2$ then 
these are the irregular classes of 
$V_1^\vee$, $\End(V_1)$, $V_1\otimes V_2$, 
$\Hom(V_1,V_2)$ respectively.
For example $\<q\>^\vee = \<-q\>$ 
and $\<q_1\>\otimes \<q_2\>$ may be computed by 
writing 
$q_1=\sum a_i t^i, q_2=\sum b_it^i$ where 
$t^r=z^{-1}$ for some integer $r\ge 1$
and then considering the Galois closed list
$q_1(\ze^a t) + q_2(\ze^b t)$ of sums of the 
various Galois conjugates, where 
$\ze=\exp(2\pi i/r)$.

The {\em levels} of a class $I$ are the nonzero slopes 
of the component circles of $\End(I)$.

\subsection{Irregular curves/wild Riemann surfaces}
A rank $n$ (bare) 
irregular curve is a triple $\bSi=(\Si,\ba,\Th)$
where $\Si$ is a smooth compact complex  
algebraic curve,
$\ba\subset \Si$ is a finite set
and $\Th$ is a rank $n$ irregular class (for each point of $\ba$).
$\bSi$ is {\em tame} if $\Th$ is tame (i.e. only involves the tame circle with multiplicity $n$ at each marked point). Thus, 
specifying a rank $n$ tame curve is the same 
as choosing a curve with marked points.

If $\Si^\circ = \Si\setminus \ba$ then any algebraic connection $(E,\nabla)\to \Si^\circ$
determines an irregular curve $(\Si,\ba,\Th)$
taking the irregular classes at each marked point.
Similarly any meromorphic connection
on a vector bundle on $\Si$ determines an irregular curve, taking its polar divisor and irregular classes.

An {\em irregular type} is similar to an irregular class but involves an ordering of the active exponents (the component circles in $I$).
Thus an irregular type determines an irregular class by forgetting the ordering.
One can then define non-bare (dressed) 
irregular curves, involving irregular types, and an ordering of the points $\ba$ 
(cf. \cite{gbs} Rmk 10.6, \cite{twcv} \S4).

\subsection{Stokes directions and exponential dominance orderings}

By looking at the exponential growth rates, there is  a partial ordering $<_d$ (exponential dominance) 
on each fibre of $\cI$, as follows. 
Suppose $d\in \partial$ and 
$q_1,q_2\in\cI_d$ are distinct then (by definition) 
$$q_1 <_d q_2$$ 
if $\exp(q_1-q_2)$ is flat (has zero asymptotic expansion) on some open sectorial neighbourhood of $d$.
As usual $i\le_d j$ means $i<_d j$ or $i=j$.
Given an irregular class $\Th$ with active exponents $I$
then $<_d$ restricts to a partial order on the fibre
$I_d$ of $I$.
Since $I$ is a finite cover, this is actually a total
order for all but a finite number of points 
$\IS\subset \partial$, the {\em Stokes directions} (or oscillating directions) of the class $\Th$.
Thus if $d\in \partial\setminus\IS$ then $I_d$ is totally ordered by $<_d$.

\subsection{Singular directions and Stokes arrows}

The points of maximal decay form a subset $\apple\subset \cI$
consisting of the points where the functions $e^q$ have maximal decay, as $q$ moves in $\<q\>$. (Sometimes they will be called p.o.m.s or {\em apples}.)
Each circle $\<q\>$ has a finite number of 
apples, except the tame circle which has none.
The {\em Stokes arrows} are the pairs
$(q_1,q_2)\in \cI\times \cI$ such that 
$\pi(q_1)=\pi(q_2)$ (so they are both in the same fibre of $\cI$) 
and $q_1-q_2\in \cI$ is a point of maximal decay.
In this case write $q_1\prec_d q_2$, where $d=\pi(q_1)$. It is viewed as an arrow from $q_2$ to $q_1$.
This defines a partial order on each fibre 
$\cI_d$ and exponential dominance refines it (if  
$q_1\prec_d q_2$ then $q_1<_d q_2$).

Given an irregular class $\Th$ with active exponents $I$
then there are only a finite number of Stokes arrows in 
$I\times I$. 
They correspond to 
the points of maximal decay of the class $\End(I)$. The Stokes arrows lie over a finite set 
$\IA\subset \partial$, the {\em singular directions} 
(or anti-Stokes directions).
The corresponding {\em Stokes quiver} at $d\in \IA$ 
is the quiver with nodes $I_d$ and arrows $\prec_d$.
As in \eqref{eq: st gp} the Stokes arrows
then determine the Stokes groups
$$\ISto_d = \ISto_{\prec_d}\subset \GL(V^0_d)$$
for any  $I$-graded local system $V^0\to \partial$.
The following lemma will be useful later.

\begin{lem}\label{lem: tord exts pord}
Suppose 
$d\in \IA$ and $i,j\in I_d$ are such that $i\prec_d j$.
Then there is an open subset $U\subset \partial$ with $d\in U$ such that $i <_ej$ for all $e\in U$. 
\end{lem}
\pf This follows as $d$ is a 
point of maximal decay for the continuous function $e^{q_i-q_j}$.
\epf
\subsection{Simple example}  \label{ssn: weber} 
Consider Weber's equation $y'' =(x^2/4 +\la)y$ 
where $\la\in \IC$
(the equation for the parabolic cylinder functions)
and the corresponding connection 
$\nabla=d-A, A=\bsmx 0 & 1 \\ x^2/4 +\la & 0 \esmx dx$.
This has just one singularity, at $x=\infty$. 
Thus $\Si=\IP^1, \Si^\circ=\IC$ and $\wh \Si$ is a closed disk with boundary circle $\partial$ (the radial compactification of $\IC$).  
\begin{wrapfigure}{l}{0.37\textwidth}
\centering
\begin{tikzpicture}[scale=0.5, level/.style={thick}]
    \draw[line width=0.2mm,color=black,samples=200,smooth,
    domain=0:2*pi
        ] plot (xy polar cs:angle=\x r,
  radius= {0.5*(6+cos((2)*(\x r)))});
  \draw[line width=0.2mm,color=black,samples=200,smooth,
    domain=0:2*pi
        ] plot (xy polar cs:angle=\x r,
  radius= {0.5*(6+cos((2)*(\x r+90)))});
    \draw[loosely dashed, line width=0.1mm,color=black,domain=0:4*pi,samples=200,smooth] plot (xy polar cs:angle=\x r,radius=0.5*6);
  \draw [ultra thick] (2.5cm,0) circle [radius=0.06];
  \draw [ultra thick] (0,2.5cm) circle [radius=0.06];
  \draw [ultra thick] (-2.5cm,0) circle [radius=0.06];
  \draw [ultra thick] (0,-2.5cm) circle [radius=0.06];
    \node[below] at (0,-4) {Stokes diagram of}; 
    \node[below] at (0,-5) {the Weber equation};
    \end{tikzpicture}
\end{wrapfigure}
A short computation, or a glance at \cite{abram-stegun} \S19.8, shows
the formal solutions at $\infty$ involve the 
multivalued functions $f_\pm=\exp(q_\pm)x^{\pm \la -1/2}$ where 
$q_\pm=\pm x^2/4$. 
The exponential factors $\exp(q_\pm)$
here 
are the main contributors to the
behaviour of solutions near $x=\infty$, and their dominance is encoded in the Stokes diagram in the figure.
(Such a diagram for the Airy equation appeared in Stokes' original paper \cite{stokes1857} and was reproduced in \cite{twcv}.)
From this we see immediately the oscillating directions 
$\IS\subset \partial$
are the four directions with argument $\pi/4 + k\pi/2$ (where
the dominance changes), and the singular directions  $\IA\subset \partial$
are the real and imaginary axes (where the ratio of dominances is
largest).  
In general such diagrams are difficult to define/draw precisely 
(especially in the multi-level case), but we can define the circles that appear, as a finite cover $\pi:I\to \partial$, 
and then view the Stokes diagram 
as a (non-intrinsic) projection of $I$ to the plane. 
In this example  
$I=\<q_+\>\sqcup \<q_-\>\to \partial$, with
each circle $\<q_\pm\>$ a trivial degree one cover.
The apples (points of maximal decay) are the four points of $I$ that project to the four marked points on the diagram.
They lie over the singular directions $\IA\subset \partial$. 
There are four Stokes arrows, 
one  over each point $d\in \IA$, from the point of maximal growth to the point of maximal decay in 
 $I_d=\pi^{-1}(d)$.
 (Fig. \ref{fig: atomic} arises for $I=\bigcup_1^3\<\al_ix^2\>$, 
 from the generic reading of the triangle $\wh A_2$,
 related to Painlev\'e 4 \cite{rsode}.)
\subsection{Ramis exponential tori}
Another topological object 
attached to an irregular class 
$\Th$ is the local system $\IT\to \partial$ of Ramis exponential tori, defined as follows.
The character lattice $X^*(\IT)\subset \cI\to \partial$
of $\IT$ is the local system of finite rank lattices (free $\IZ$-modules) generated by the active exponents 
$I\subset \cI$, so that 
$X^*(\IT)_d = \< I_d \>_\IZ\subset \cI_d$.
Then $\IT$ is the local system of 
tori with this character lattice, so that 
$\IT_d = \Hom(X^*(\IT)_d,\IC^*)$.
Note that if $V^0\to \partial$ is an $I$-graded 
local system of vector spaces with dimension $\Th$, then
$V^0$ can be viewed as graded by $X^*(\IT)$.
This means there is a faithful action of $\IT$ on $V^0$, i.e. an injective map $\IT\to \GL(V^0)$ of local systems of groups.
(These tori appear in the differential Galois group of the corresponding connections.)
\section{Stokes filtered local systems}
\label{sn: SFLS}

Recall that an irregular class $\Th:\cI\to \IN$ determines the data 
$n,I, \IS, <_d$.

A {\em Stokes filtration} $F$ of type $\Th$ on a local system 
$V\to  \wh\Si$ of complex vector spaces
is the data of an $I_d$-filtration 
$F_d$ of $V_d$ of dimension $\Th$, 
for each $d\in \partial\setminus\IS$, such that:

{\bf SF1)} The $F_d$ are locally constant as $d$ varies in $\partial$ 
without crossing $\IS$,

{\bf SF2)} The Stokes condition \eqref{eq: Fstokescondns} holds across each $d\in \IS$.

To be precise, in 2)  the filtrations on the left and right of $d\in \IS$ are transported 
to the fibre $V_d$ to get into the exact situation of \eqref{eq: Fstokescondns}.
Thus in other words 
condition 2) says that there 
exists a local grading inducing the filtrations
(i.e. a local splitting across $d$): 
there is a grading $\Ga$ of 
the local system $V$ (by sub-local systems) 
throughout a small neighbourhood 
$U$ of $d$ in $\partial$, such that  
$F_e=\cF(\Ga,<_e)$ for all $e\in U\setminus \{d\}$.
The condition 1) means parallel transport along any path in $\partial\setminus\IS$ relates the filtrations.

\begin{defn}
A Stokes filtered local system is a triple $(V,\Th, F)$
where $V\to \wh\Si$ is a local system of vector spaces,
$\Th$ is an irregular class (of the same rank as $V$)
and $F$ is a Stokes filtration on $V$ of type $\Th$.
\end{defn}

\begin{rmk}[Robustness]\label{rmk: robustness}
Note that the definition is robust in the sense 
that if a finite number of points is added to $\IS$ 
(and the filtrations $F_d$ 
are not specified at these points), 
then nothing is changed since SF2) implies $F$ is
continuous across the missing points, and thus determined by neighbouring filtrations.
\end{rmk}

Two Stokes filtered local systems 
$\cV_i=(V_i,\Th_i, F^i)$ for $i=1,2$ are isomorphic if $\Th_1=\Th_2$ and there is an isomorphism
$\varphi : V_1\to V_2$ of local systems 
relating the Stokes filtrations.

To define more general morphisms 
first note that if 
$(V,\Th, F)$ is a Stokes filtered local system and 
$d\in \partial\setminus \IS$ and $i\in \cI_d$ then 
one can define
\beq\label{eq: full stokes filtration}
F_d(i) = \sum_{j\in I_d\ \vert\  j\le_d i} F_d(j) \subset V_d.
\eeq
Clearly if $i\le_d j$ then $F_d(i)\subset F_d(j)$.
Thus $V_d$ can be viewed as filtered by the 
poset $\cI_d$ and not just by the ordered set $I_d$.
A morphism
$\cV_1\to \cV_2$
of Stokes filtered local systems
is a morphism 
$\varphi : V_1\to V_2$
of local systems
that restricts to a map of $\cI_d$-filtered vector spaces 
($\varphi(F^1_d(i))\subset F^2_d(i)$ for all $i\in \cI_d$),  
for all $d\in \partial\setminus(\IS_1\cup \IS_2$).

If $\cV=(V,\Th,F)$ is a Stokes filtered local system, define a global section of $\cV$
to be 
$$v\in \HH^0(\cV,\wh\Si) := \{ v\in \HH^0(V,\wh \Si)\st v(d)\in 
F_d\<0\>
\text{ for all } d\in \partial\setminus \IS\}$$
i.e. it is a global section of $V$
such that $v(d)\in F_d\<0\>:=F_d(\<0\>)$ for all points 
$d\in \partial\setminus \IS$.  
Said differently, using the Stokes filtrations, 
 define a section of 
$V$ to have ``moderate growth'' in the direction 
$d\in \partial$ if $v(e)\in F_e\<0\>$
for all non-Stokes direction $e$ in some open neighbourhood of $d$. 
Then $\HH^0(\cV,\wh\Si)$ is just the space of 
sections of $V$ that have moderate growth everywhere.
It will become clear later (Prop. \ref{prop: fflem})
that $\cV$ has no global sections unless the tame circle $\<0\>$ is active in the irregular class at each marked point.

The dual Stokes filtered local system
$\cV^\vee$ is defined as 
$(V^\vee, \Th^\vee, F^\vee)$ where 
 $$F_d^\vee(q) = F_d(<\!\!-q)^\perp\subset V_d^\vee,$$
and the tensor product $\cV_1\otimes \cV_2$
is $(V_1\otimes V_2, \Th_1\otimes \Th_2, F)$
where 
$$F_d(q) = \sum_{q=q_1+q_2} F^1_d(q_1)\otimes F^2_d(q_2)$$
for generic $d$ (but this is sufficient by Rmk \ref{rmk: robustness}).
Thus  one can define 
$\Hom(\cV_1,\cV_2)=\cV_2\otimes \cV_1^\vee$
and observe that
a morphism $\cV_1\to \cV_2$ 
is the same thing as a global section of this Stokes filtered local system.

The definition used by Deligne and Malgrange is easily seen to be equivalent to this
(this amounts to recognising, as in \S\ref{ssn: int filtr} below, that filtrations of $V_d$
by $\cI_d$ for $d\in \IS$ are canonically 
determined too). Thus (the global version of) their Riemann--Hilbert--Birkhoff correspondence says that:

\begin{thm}\label{thm: delmalg}
The category of Stokes filtered local systems is equivalent
to the category of algebraic connections $(E,\nabla)$
on vector bundles on $\Si^\circ$.
\end{thm}

This was 
conjectured by Deligne \cite{deligne78in}
and proved  by Malgrange \cite{Mal-irregENSlong} 
(using the earlier Malgrange--Sibuya cohomological classification).
As remarked in \cite{malg-book} p.57, this is essentially equivalent to the 
statement of Jurkat \cite{jurkat78}---this equivalence will be discussed  in the sequel.
The definition of the Stokes filtration of a connection
will be discussed in Apx. \ref{sn: anbbx}.

\subsection{} \label{sn: weberfilt}
\begin{wrapfigure}{r}{0.37\textwidth}
\centering
\begin{tikzpicture}[scale=0.5, level/.style={thick}]
    \draw [thick, white] (0,0) circle [radius=4.2cm];
    \draw [thick] (0,0) circle [radius=4cm];
  \draw [ultra thick] (2.828cm,2.828) circle [radius=0.06];
  \draw [ultra thick] (2.828cm,-2.828) circle [radius=0.06];
  \draw [ultra thick] (-2.828cm,2.828) circle [radius=0.06];
  \draw [ultra thick] (-2.828cm,-2.828) circle [radius=0.06];
  \node[left] at (4,0) {$L_0$};
  \node[right] at (-4,0) {$L_2$};
  \node[below] at (0,4) {$L_1$};
  \node[above] at (0,-4) {$L_3$};
  \node[left] at (-3.7,1.4) {$\partial$};
  \node[below] at (0,-4) {$\wh \Si$ with $\IS\subset \partial$ marked};

    \end{tikzpicture}
\end{wrapfigure}
Returning to the simple example of \S\ref{ssn: weber}, 
the Stokes filtrations (in this example) 
amount to recording the lines in $V$ spanned by the
subdominant (or recessive) solutions on each sector 
$\partial\setminus \IS$.
Intrinsically the filtrations are indexed by the components of 
$I\setminus\pi^{-1}(\IS)$ (totally ordered over each sector by the dominance of the exponential factors).
Since $V$ is trivial this amounts to specifying four one dimensional subspaces
$L_0,L_1,L_2,L_3\subset \HH^0(V)\cong \IC^2$
where $L_k$ is the line of
solutions which are recessive at $\infty$ when $\arg(x)=k\pi/2$.
The Stokes condition \eqref{eq: Fstokescondns} means that $L_i\neq L_{i+1}$ for all $i$ (indices modulo $4$).

\section{Stokes graded local systems}\label{sn: SGLS}

Recall that an irregular class $\Th:\cI\to \IN$ determines the data 
$n,I, \IA, \prec_d$.

A {\em Stokes grading} $\Ga$ of type $\Th$ of a 
local system $V\to  \wh\Si$ 
is the data of an $I_d$-grading 
$\Ga_d$ of $V_d$ of dimension $\Th$, 
for each $d\in \partial\setminus\IA$, such that:

{\bf SG1)} The $\Ga_d$ are locally constant as $d$ varies in $\partial$ 
without crossing $\IA$,

{\bf SG2)} The Stokes condition \eqref{eq: Gstokescondns} 
holds across each $d\in \IA$.

To be precise, in 2) the gradings on the left and right of $d\in \IA$ are transported to the fibre $V_d$ 
to get into the exact situation of \eqref{eq: Gstokescondns}. 
 
\begin{defn}
A Stokes graded local system is a triple $(V,\Th, \Ga)$
where $V\to \wh\Si$ is a local system of vector spaces,
$\Th$ is an irregular class (of the same rank as $V$)
and $\Ga$ is a Stokes grading on $V$ of type $\Th$.
\end{defn}

Again the definition is robust, by Rmk \ref{rmk: empty prec}.

Note that
if $d\in \IA$ and $L,R\in \partial$ are points just to the left and right of $d$, then there are {\em two} distinguished isomorphisms
$V_L\to V_R$: one given by the transport of the local system $V$, and the other
given by the graded isomorphism ({\em wild transport})
\beq\label{eq: tspt across singdir}
(V_L,\Ga_L)\cong 
(V_d,\Ga_L)\mapright{g}
(V_d,\Ga_R) \cong (V_R,\Ga_R)
\eeq
where the first and third isomorphisms come from the local system
structure of $V$, and $g=g(\Ga_L,\Ga_R)$ is the wild monodromy. This leads to  the
Stokes local system (\S\ref{sn: SLS}).

The terms {\em Stokes grading} or {\em Stokes decomposition} will be used interchangeably.
If $d\in \IA$ then the fibre $V_d$ can be given the median grading of the fibres to either side
 (cf. \S\ref{sn: sls and sgls}), yielding a 
 preferred decomposition of each tangential fibre of $V$.

The global sections of a 
Stokes graded local system
$\cV=(V,\Th,\Ga)$
are the sections of the underlying local system $V$
that go into the piece graded by the tame circle
$\<0\>$ in 
each singular sector:
\beq\label{eq: sgls sns}
\HH^0(\cV,\wh \Si) := 
\{ v\in \HH^0(V,\wh \Si)\st v(d)\in \Ga_d\<0\>
\text{ for all } d\in \partial\setminus \IA\}.
\eeq

A morphism $\cV_1\to \cV_2$ 
is a morphism 
$\varphi : V_1\to V_2$
of local systems
that restricts to a map of $\cI_d$-graded 
vector spaces 
($\varphi(\Ga^1_d(i))\subset \Ga^2_d(i)$ 
for all $i\in \cI_d$),  
for all $d\in \partial\setminus(\IA_1\cup \IA_2$).

The dual Stokes graded local system
$\cV^\vee$ is defined as 
$(V^\vee, \Th^\vee, \Ga^\vee)$ where 
 $$\Ga_d^\vee(q) = 
 \left(\bigoplus_{i\in I_d\setminus\{-q\}} \Ga_d(i)\right)^\perp\subset V_d^\vee,$$
and the tensor product $\cV_1\otimes \cV_2$
is $(V_1\otimes V_2, \Th_1\otimes \Th_2, \Ga)$
where 
$$\Ga_d(q) = \bigoplus_{q=q_1+q_2} 
\Ga^1_d(q_1)\otimes \Ga^2_d(q_2)$$
for generic $d$ (but this is sufficient by robustness).
Thus  one can define 
$\Hom(\cV_1,\cV_2)=\cV_2\otimes \cV_1^\vee$
and observe that
a morphism $\cV_1\to \cV_2$ 
is the same thing as a global section of this Stokes graded local system.

\subsection{}  \label{ssn: weberstgr}

\begin{wrapfigure}{r}{0.37\textwidth}
\centering
\begin{tikzpicture}[scale=0.5, level/.style={thick}]
    \draw [thick, white] (0,0) circle [radius=4.2cm];
    \draw [thick] (0,0) circle [radius=4cm];
  \draw [ultra thick] (0,4) circle [radius=0.06];
  \draw [ultra thick] (0,-4) circle [radius=0.06];
  \draw [ultra thick] (4,0) circle [radius=0.06];
  \draw [ultra thick] (-4,0) circle [radius=0.06];
  \node[left] at (3,2.4) {$\Ga_1$};
  \node[right] at (-2.9,2.4) {$\Ga_2$};
  \node[right] at (-3,-2.4) {$\Ga_3$};
  \node[left] at (3,-2.4) {$\Ga_4$};
  \node[left] at (-3.7,1.4) {$\partial$};
  \node[below] at (0,-4) {$\wh \Si$ with $\IA\subset \partial$ marked};
    \end{tikzpicture}
\end{wrapfigure}
Returning to the simple example of 
\S\ref{ssn: weber}, 
a Stokes grading of type $I=\<q_-\>\cup\<q_+\>$ 
amounts to specifying a direct sum decomposition
of $V$ on  each sector of $\partial\setminus \IA$.
Intrinsically the Stokes gradings are indexed by the components of 
$I\setminus\pi^{-1}(\IA)$.
Thus (away from $\IS$) this amounts to choosing
two complementary one dimensional subspaces of $V$, 
that we call the ``subdominant'' and ``dominant'' solutions. 
In this example the
Stokes conditions \eqref{eq: Gstokescondns} 
say that at each point $d\in \IA$
the lines of ``subdominant'' solutions on each side match up.
Consequently any ``dominant'' solution may jump at $d$, by the addition of a ``subdominant'' solution.

\section{Stokes local systems}\label{sn: SLS}

As in \S\ref{sn: SGLS}, 
recall that an irregular class $\Th:\cI\to \IN$ determines the data 
$n,I, \IA, \prec_d$. 
Also for
$d\in \IA$ the Stokes arrows $\prec_d$ 
determine the Stokes group 
$\ISto_d\subset \GL(V_d^0)$ for any $I_d$-graded vector space $V_d^0$, as in \eqref{eq: st gp}.

Two slightly different definitions of a Stokes local system will be given, the first with explicit gluing maps.

Suppose $V\to  \wh\Si$ is a rank $n$ 
local system of complex vector spaces and 
$V^0\to \partial$ is an $I$-graded local system of dimension $\Th$.
Thus on $\partial$ there are two local systems, $V^0$ and the restriction of $V$.

A {\em gluing map} $\Phi$ is an isomorphism of vector spaces 
$\Phi_d:V^0_d\cong V_d$ for each $d\in \partial\setminus \IA$, which is locally constant as $d$ varies in $\partial$ without crossing $\IA$.
Thus it is a section of the local system 
$\Iso(V^0,V)\subset \Hom(V^0,V)$ over 
$\partial\setminus \IA$.

A gluing map determines an automorphism 
$\IS_d\in \GL(V_d^0)$ for each $d\in \IA$, as the composition of the isomorphisms:
\beq\label{eq: first StoAM}
V_d^0\ \mapright{1} \ 
V_R^0 \ \mapright{\Phi_R}\  
V_R \ \mapright{3}\  
V_L \ \mapright{\Phi_L^{-1}}\  
V_L^0 \ \mapright{5} \  V_d^0
\eeq
where $L/R$ denotes fibres just to the left/right of $d$, the maps $1,5$ are the transport of $V^0$ and 3 is the transport of $V$. (More precisely $L$ means the positive side of $d$).

A gluing map is {\em Stokes} if $\IS_d\in\ISto_d$ for all $d\in \IA$.

\begin{defn}
A ``Stokes local system with gluing maps'' 
is a $4$-tuple $(V,V^0,\Th, \Phi)$
where $V\to \wh\Si$ is a local system of vector spaces,
$\Th$ is an irregular class (of the same rank as $V$),
$V^0\to \partial$ is an $I$-graded local system of dimension $\Th$,
and $\Phi$ is a Stokes gluing map.
\end{defn}

Note that the gluing map 
$\Phi$ is Stokes if and only if $(V,\Th, \Ga)$
is a Stokes graded local system, where
$\Ga$ is the set of gradings 
defined by $\Ga_d(i)=\Phi_d(V^0_d(i))$
for all $i\in I_d, d\in \partial\setminus \IA$.
This defines a functor to Stokes graded local systems,
and it will be seen below that it is an equivalence.

The wild monodromy at $d\in \IA$ of a Stokes
local system with gluing maps is the
automorphism 
$g_d\in \GL(V_d)$ given by 
the composition of the isomorphisms:
\beq\label{eq: wildmonod}
V_d\ \mapright{1} \ 
V_L \ \mapright{\Phi_L^{-1}}\  
V_L^0 \ \mapright{3}\  
V_R^0 \ \mapright{\Phi_R}\  
V_R \ \mapright{5} \  V_d
\eeq
where $L/R$ denotes fibres just to the left/right of $d$, the maps $1,5$ are the transport of $V$ and 3 is the transport of $V^0$.
Clearly $g_d$ is conjugate to $\IS_d$, but 
acts on a different space.

By definition 
the global sections of $\cV=(V,V^0,\Phi,\Th)$
are pairs $(v,v^0)$ where $v\in \HH^0(V,\wh \Si)$
and $v_0\in \HH^0(V^0\<0\>,\partial)$,
which map to each other under the gluing maps (at each point $d\in \partial\setminus \IA$).
This is the same as giving just 
$v\in \HH^0(V,\wh \Si)$ such that 
$v(d)$ is fixed by $g_d$ for all $d\in \IA$, and 
$\Phi_d^{-1}(v(d))\in V_d^0\<0\>$ for all 
$d\in \partial\setminus \IA$.
However it follows from 
Lemma \ref{lem: monod invt vectors} that the first condition is redundant and so 
\beq\label{eq: sls sns}
\HH^0(\cV,\wh \Si) = 
\{ v\in \HH^0(V,\wh \Si)\st v(d)\in \Phi_d(V^0_d\<0\>)
\text{ for all } d\in \partial\setminus \IA\}.
\eeq

To motivate the second approach,
note that  the graded local system $V^0\to \partial$ is the same as 
giving a graded local system on each halo (a germ of a punctured disk around each puncture). 
Thus there are two local systems on each halo: $V^0$ and the restriction of $V$. The gluing maps glue them together, away from each singular direction.
One can picture this by putting $V^0$
on a second copy of the halo, above the first one 
(with $V$ on it), and then gluing the two copies, away from the 
singular directions (leading to a small tunnel over each 
singular direction). 
The second version appears by removing the base of each tunnel, and then pushing this picture flat (yielding a surface with tangential punctures).

Said differently for the second version just 
restrict $V$ to the complement of the halos,
so the gluing maps amount to a way to glue $V^0$ to $V$ across each component of the outer boundary of each halo, away from the singular directions.
This is the most convenient description in practice.
It leads to the idea of  {\em fission} \cite{fission}, 
that the structure group is broken 
(to the group of graded automorphisms) 
at the boundary of $\wh \Si$, as illustrated in the picture on the title page. (Here the $I$-graded local system is identified with a local system on $I\to \IH$  in the usual way, and the cover $I$ is glued to the interior of the curve---in other words a Stokes local system is thus the same as a 
sort of generalised local system on the surface in the picture, where the ranks jump as the surface bifurcates. The tangential punctures are not drawn but should be there). 

Define the halos $\IH \subset \wh\Si$ 
as a small tubular neighbourhood of $\partial$, so 
$\IH$ is a union of annuli. 
Let $\partial'$ be the boundary of $\IH$ in the interior of $\wh \Si$. 
The cover $I\to \partial$ extends to a cover $I\to \IH$.
Choose a smooth bijection $e:\partial\to \partial'$ preserving  the order of the points.
Thus if the irregular class $\Th$ and thus  
$\IA\subset \partial$ is given, 
then  the auxiliary surface:
$$\wt \Si = \wt\Si(\Th) := \wh \Si\setminus 
e(\IA)$$
is well defined, 
removing the tangential punctures
$e(d)$ from $\wh\Si$ for $d\in \IA$.
If $d\in \IA$ let $\ga_d$ be the small positive loop in 
$\wt \Si$ based at $d$ 
that goes across $\IH$, around $e(d)$ and back to $d$.

\begin{defn}
A Stokes local system 
is a pair $(\IV,\Th)$
where 
$\Th$ is an irregular class and
$\IV$ is a local system of vector spaces on
$\wt \Si(\Th)$, 
equipped with an $I$-grading over $\IH$
(of dimension $\Th$),
such that the monodromy 
$\IS_d:=\rho(\ga_d)$ is in 
$\ISto_d\subset \GL(\IV_d)$ 
for each $d\in \IA$.
\end{defn}

The element $\IS_d=\rho(\ga_d)\in
\ISto_d\subset \GL(\IV_d)$ is the Stokes automorphism.
The equivalence of the two approaches is straightforward, 
whence $\IS_d$ 
is identified with that defined in \eqref{eq: first StoAM}.
This is the specialisation to general linear groups of 
the definition in \cite{gbs,twcv}.

Note that an $I$-graded local system of dimension $\Th$
is the same thing as an $\cI$-graded local system of dimension $\Th$
(all the other components of $\cI$ grade the trivial rank zero sublocal system).

If $\cV=(\IV,\Th)$ is a Stokes local system then a section of $\cV$ 
is a section of $\IV$ that takes values in $\IV\<0\>$  over each halo 
(the graded piece indexed by the tame circle). 

A morphism of Stokes local systems $\cV_1\to \cV_2$ 
is a map of local systems 
on $\wt \Si(\Th_1)\cap\wt\Si(\Th_2)$ that restricts to a map of $\cI$-graded local systems over $\IH$.
In particular two  Stokes local systems are isomorphic 
if and only if their irregular classes are 
equal and their
underlying local systems $\IV_1,\IV_2$ are isomorphic.

Note on terminology:
The ``wild monodromy'' and 
``Stokes automorphism''
are essentially the same thing---namely 
the monodromy of the 
Stokes local system around a small positive loop 
around a tangential puncture.
In this paper the term 
``Stokes automorphism'' will be used 
if the loop is based in the halo, and 
the term ``wild monodromy'' will be used 
if the loop is based outside the halo (in the interior of the curve). 
For example $g_d\in \GL(V_d)$ 
and $\IS_d\in \GL(V^0_d)$.
Note however that these {\em equal}
if one starts with a Stokes graded local system
(see Lemma \ref{eq: wm and sam}).

\subsection{} \label{ssn: webersls}
\begin{wrapfigure}{r}{0.37\textwidth}
\centering
    \begin{tikzpicture}[scale=0.5, level/.style={thick}]
    \draw [fill={rgb:black,1;white,8}, thick] (0,0) circle [radius=4cm];
    \draw [fill=white, thick, color=white] (0,0) circle [radius=2cm];
    \draw[thick,color=gray,samples=200,smooth,
    domain=0:2*pi
        ] plot (xy polar cs:angle=\x r,
  radius= {0.5*(6+cos((2)*(\x r)))});
  \draw[thick, color=gray,samples=200,smooth,
    domain=0:2*pi
        ] plot (xy polar cs:angle=\x r,
  radius= {0.5*(6+cos((2)*(\x r+90)))});
  \draw [thick,fill=white] (2cm,0) circle [radius=0.15];
  \draw [thick,fill=white] (0,2cm) circle [radius=0.15];
  \draw [thick,fill=white] (-2cm,0) circle [radius=0.15];
  \draw [thick,fill=white] (0,-2cm) circle [radius=0.15]; 
  \draw[->, ultra thick, gray] (3.5cm,0cm) -- (2.5cm,0cm);%
  \draw[->, ultra thick, gray] (-3.5cm,0cm) -- (-2.5cm,0cm);
  \draw[->, ultra thick, gray] (0cm,3.5cm) -- (0cm,2.5cm);
  \draw[->, ultra thick, gray] (0cm,-3.5cm) -- (0cm,-2.5cm);
  \node[left] at (-2.9,2.9) {$\partial$};
  \node[below] at (0,-4) {$\wt \Si$ with Stokes diagram};
  \node[below] at (0,-5.2) {and Stokes arrows drawn};
    \end{tikzpicture}
\end{wrapfigure}
Returning to the simple example of 
\S\ref{ssn: weber}, 
a Stokes local system of type $I=\<q_-\>\cup\<q_+\>$ 
is a local system $\IV$ on the auxiliary surface $\wt \Si$,
graded by $I$ over the halo (the shaded area in the figure).
Given $d\in \IA$ the fibre $\IV_d$ is graded by $I_d$
and so the Stokes arrow between the two points of 
$I_d$ determines the Stokes group
$\bsmx 1 & * \\ 0 & 1 \esmx \subset \GL(\IV_d)$, 
where the star corresponds to maps along the arrow. 
The monodromy of $\IV$ around a small loop (based at $d$)
around the tangential puncture $e(d)$,  
should be in this Stokes group.
Choosing a basepoint $b\in \partial$ and suitable loops in 
$\wt \Si$ based at $b$ then yields the formal monodromy
$h\in \GrAut(\IV_b)$ (the monodromy of $\IV$ around $\partial$)
and Stokes automorphisms $S_1,S_2,S_3,S_4$ (in alternating unipotent groups) satisfying the monodromy relation  $hS_4S_3S_2S_1=1$.

\begin{rmk}
Here are more details on how to see the two versions are equivalent.
Choose a retraction of $\IH\cong [0,1]\times \partial$
onto $\partial$  (dragging 
$e(d)$ along a path $\la_d$ to $d$).
Consider the intermediate 
category of four-tuples $(V,V^0,\Phi,\Th)$ where $V^0$ is now defined on all of $\IH$, and the gluing maps are 
defined on all of $\IH$ except on the paths (cilia) 
$\la_d$ for $d\in \IA$.
Such  objects restrict to give a Stokes local system with gluing maps in the obvious way, yielding an equivalence.
On the other hand the restriction of $V$ to the complement of $\IH$ really does now glue to $V^0$ on 
$\partial' \setminus e(\IA)$,  
defining a Stokes local system $\IV$, 
yielding an equivalence.
\end{rmk}

\begin{rmk}
1) An alternative version (cf. \cite{p12} Apx. B) 
is to glue the halo 
$\IH=\partial\times [0,1]$ on the {\em other side}  of $\partial$, and then the tangential punctures are made at 
$\IA\subset \partial$. 
This is homotopy equivalent to the presentation above.
In practice it is not important 
exactly where the tangential punctures are, 
provided they are in the right order 
(cf. also the version in \cite{gbs} with non-crossing cilia drawn to keep them in order, which can be pulled tight whenever convenient---this perhaps best reflects the desired picture especially when the curve etc moves).

2) Perhaps the simplest approach is to just glue a second copy of 
$\partial$ onto $\partial$ away from the singular directions 
$\IA\subset \partial$. One can prove the resulting (non-Hausdorff) space has the desired fundamental group \cite{MO-hatcher}, and then view $\IV$ as a local system on it.
\end{rmk}

\section{Stokes local systems and Stokes
graded local systems}\label{sn: sls and sgls}

A Stokes local system with gluing maps 
$(V,V^0,\Phi,\Th)$ determines a Stokes grading
$\Ga$,  by taking the 
image of the graded local system under the gluing maps: 
$$\Ga_d(i)=\Phi_d(V^0_d(i))$$  
for all $d\in \partial\setminus \IA, i\in I_d$.
Conversely given a Stokes graded local system
$(V,\Th, \Ga)$
it is easy to construct a 
Stokes local system (with gluing maps)
mapping to it as above:
If $d\in \IA$ then transport the gradings on either side of $d$ to $V_d$ and define $\Ga_d$ to be the
median grading (of $V_d$) determined by these gradings from either side, as in \S\ref{sn: median gr}.
Then, for any $d\in \partial$, define
$V^0_d$ to be the graded vector space $(V_d,\Ga_d)$.
The spaces $V^0_d$ 
form the fibres of a graded local system $V^0\to \partial$: for any path in $\partial\setminus \IA$ this is clear. 
On the other hand for a small path 
across some $d\in \IA$ use 
the ``wild transport'' \eqref{eq: tspt across singdir}.
Similarly for a path ending at $d\in \IA$, 
for example:
$$(V_L,\Ga_L)\ \cong\  
(V_d,\Ga_L)\ \mapright{g(\Ga_L,\Ga_d)}\ 
(V_d,\Ga_d).$$

This defines the graded local system $V^0$, and by construction it comes with gluing maps (the identity)
$\Phi_d:V^0_d\to V_d$ for all 
$d\in \partial\setminus  \IA$.

To show this is an equivalence of categories it remains to check fully faithfulness, which is left as an exercise.
(Using the internal hom it is enough to show the map between spaces of global sections is bijective, and 
this is immediate, comparing 
\eqref{eq: sgls sns} and  \eqref{eq: sls sns}.)

\begin{lem}\label{eq: wm and sam}
Suppose $(V,V^0,\Phi,\Th)$ is the Stokes local system
with gluing maps 
of a Stokes graded local system $(V,\Ga,\Th)$
as constructed above.
Then the wild monodromy equals the Stokes automorphism
as an element of $\GL(V_d)$ for any $d\in \IA$.
\end{lem}
\pf
Recall that $V_d^0$ is the space $V_d$ equipped with the median grading $\Ga_d$, so the statement makes sense.
There are three gradings $\Ga_L,\Ga_R,\Ga_d$ in $V_d$. 
Then $g_d = g(\Ga_L,\Ga_R)$, and 
the definition of $\IS_d$ says that 
$\IS_d = g(\Ga_L,\Ga_d)\circ g(\Ga_d,\Ga_R)= (\sqrt{g_d})^2=g_d$.
\epf

\section{Operations on Stokes filtered local systems}
\subsection{Intermediate filtrations}\label{ssn: int filtr}

Suppose $(V,\Th,F)$ is a Stokes filtered local system indexed by $I$.
If $d\in \partial$ is a Stokes direction for $I$
then the dominance ordering of $I_d$ is only a partial order.
Nonetheless for each $i\in I_d$ a 
subspace $F_d(i)\subset V_d$ can be defined as follows.
Let $F_L(i),F_R(i)$ be the corresponding steps of the 
Stokes filtrations, just on the left and right of $d$ respectively.
Transport them to $d$ (using the local system structure of $V$),
to obtain subspaces 
$F_{L/R}(i)\subset V_d$, and define the {\em intermediate filtration}:
\beq\label{eq: intfiltration}
F_d(i) := F_L(i)\cap F_R(i).
\eeq
These subspaces have the property that if $i<_d j$ then
$F_d(i)\subset F_d(j)$. 
Note that the partial order $<_d$ is the intersection of the adjacent orders on either side of $d$.

\begin{lem}\label{lem: intfilt}
If $d\in \IS, i\in I_d$ and 
$\Ga$ is any local splitting of $F$ across $d$, then
$F_d(i)$ equals the filtration associated to the grading $\Ga$ by the partial order $<_d$:
$$F_d(i) = \cF(\Ga,<_d)(i):=\Ga_d(i)\oplus\bigoplus_{j<_d i} \Ga_d(j).$$
\end{lem}
\pf
Fix $d\in \IS$ and $i\in I_d$.
Then $i$ determines a 
partition of $I_d$ into five subsets: 
$$I_d=I_{++}\sqcup I_{--}\sqcup I_{-+}\sqcup I_{+-}\sqcup I_{=},$$ where
$I_==\{i\}$ and for example $I_{+-}$ is the subset of elements of $I_d$ that are greater than $i$ on the left and less than $i$ on the right of $d$.
Thus the choice of any local splitting  determines a decomposition
$V_d=V_{++}\oplus V_{--}\oplus V_{-+}\oplus V_{+-}\oplus V_{=}$ (summing over the corresponding indices).
By construction $F_L = V_{--}\oplus V_{-+}\oplus V_{=}$,
$F_R = V_{--}\oplus V_{+-}\oplus V_{=}$ and 
$F_d(i) = V_{--}\oplus V_{=}$, so the claim follows:
$F_d(i)=F_L(i)\cap F_R(i)$. 
\epf

Thus one can just as well incorporate the data of the 
filtration of every fibre $V_d$ (not just away from the Stokes directions) in this slightly generalised sense (since $I_d$ is only partially ordered, they are not quite filtrations in the usual sense).
As shown above this ``extra data'' is canonically determined by the Stokes filtration.

\subsection{Associated graded local system}

Suppose $(V,\Th, F)$ is a Stokes filtered local system with active exponents $I$.
For each $d\in \partial$ one can consider the 
associated graded vector space $\Gr(V_d,F_d)$.
Even though the filtrations will in general jump discontinuously at
Stokes directions, these graded vector spaces fit together in a 
canonical way as the fibres of a graded local system:

\begin{lem}
Given a Stokes filtration $F$ on a local system $V\to \partial$ indexed by $I\subset \cI$, then there is a canonically determined 
$I$-graded local system, $\Gr(V,F)\to \partial$, the associated graded, with fibres the associated graded vector spaces $\Gr(V_d,F_d)$.
\end{lem}
\pf
Away from Stokes directions this is immediate 
since the Stokes filtrations are locally constant. 
Suppose $d$ is a Stokes direction.
In the notation of the proof of Lemma \ref{lem: intfilt}
the piece of $\Gr(V_d,F_d)$ with index $i\in I_d$ 
is  $F_d(i)/F_d(<\!\!i)$ at $d$, and on the left the corresponding graded piece is
$F_L(i)/F_L(<\!\!i)$.
Given a local grading then 
$F_d(<\!\!i) = V_{--}$ and 
$
F_L(<\!\!i) = V_{--}\oplus V_{-+}.$
It follows that the natural inclusion map 
$F_d(i)=F_L(i)\cap F_R(i)\into F_L(i)$ 
maps $F_d(<\!\!i)$ 
into $F_L(<\!\!i)$, and induces an isomorphism  
$F_d(i)/F_d(<\!\!i)\cong F_L(i)/F_L(<\!\!i)$ on the quotients.
Similarly going to the right of $d$.
This provides the gluing maps to define the graded local system.
\epf

In the simple example of 
\S\S\ref{ssn: weber},\ref{sn: weberfilt},
the associated graded $V^0=\Gr(V,F)\to \partial$
is isomorphic to the local system whose solutions
are the two functions 
$f_\pm=\exp(q_\pm)x^{\pm \la-1/2}$
appearing in the formal solutions at $\infty$. 
The monodromy of $V^0$ (the formal monodromy) 
is thus $-\diag(t,t^{-1})$ where $t=\exp(2\pi i \la)$.
A nice exercise shows this amounts 
to taking the cross-ratio of the four lines 
$L_1,\ldots,L_4\subset \IC^2$ (see \cite{even-euler} Lem. 11).

\section{Stokes filtrations from Stokes gradings}\label{sn: SFLS from SGLS}

A Stokes graded local system $(V,\Th,\Ga)$ determines a
Stokes filtered local system $(V,\Th,F)$
by taking the associated filtrations, 
using the dominance orderings.
In detail this goes as follows.
For any $d\in \partial\setminus \IS$ 
we need to define an $I_d$-filtration $F_d$ of $V_d$.
There are two cases:

1) If $d\not\in \IA$ then  $F_d=\cF(\Ga,<_d)$
is the filtration associated to the grading using the 
dominance ordering of $I_d$.

2) If $d\in \IA$ then transport the two gradings on 
either side, to $V_d$ and then take their associated filtrations. 
The Stokes condition ensures these two filtrations are equal,
since the order $<_d$ extends the partial order $\prec_d$,
 by  Lemma \ref{lem: tord exts pord}.

To show this is a Stokes filtration first note that 
Lemma \ref{lem: tord exts pord}
implies the following:

\begin{cor}[Malleability]\label{cor: malleability}
Suppose $(V,\Th,\Ga)\to \partial$ is a Stokes graded local system
and $d\in \IA\cap \IS$.
Let $\Ga_L,\Ga_R$ be the gradings just to the left and right of 
$d$ and let $F_L,F_R$ be the corresponding associated filtrations.
Transport all this to $V_d$.
Then both gradings split both filtrations.
\end{cor}
\pf By Lem. \ref{lem: tord exts pord} 
the orders $<_L,<_R$ from either side 
extend the partial order $\prec_d$.
Thus the result follows immediately from the Stokes conditions at $d$ for $\Ga_L,\Ga_R$.
\epf

\begin{lem}
$(V,\Th,F)$ is a Stokes filtered local system.
\end{lem}
\pf
For any $d\in \IS$ we need to  check the 
two filtrations on either side of $d$ satisfy the Stokes condition. 
Firstly if $d\not\in \IA$ then this is clear, since then 
$\Ga_d$ gives the desired splitting.
Secondly if $d\in \IA$, then the transport of 
the grading on either side will split both filtrations, by 
Corollary \ref{cor: malleability}.
\epf

This defines a functor $\varphi$ from Stokes graded local systems to Stokes filtered local systems.
It is fully faithful and surjective so is
an  equivalence of categories 
(it is surjective and not just essentially surjective).
Surjectivity follows from the main theorem:

\begin{thm} \label{thm: v2 of main thm}
Any Stokes filtered local system $(V,\Th,F)$ 
admits a unique Stokes grading $\Ga$, such that 
$F$ is the Stokes filtration associated to $\Ga$:
$$(V,\Th,F)=\varphi(V,\Th,\Ga).$$
\end{thm}

This will be proved in the next section, and fully faithfulness in \S\ref{sn: ff}.

\begin{rmk}\label{rmk: same gls}
Note that a Stokes graded local system thus determines a graded local system in two ways: as the associated graded of the associated filtrations, or via the converse construction in 
\S\ref{sn: sls and sgls}. The fact they  are canonically isomorphic 
is left as an exercise.
\end{rmk}

\subsection{}
In the simple example of 
\S\S\ref{ssn: weber},\ref{sn: weberfilt},\ref{ssn: weberstgr}, 
Thm. \ref{thm: v2 of main thm} holds since 
the Stokes gradings are then given by the two adjacent subdominant solutions: $\Ga_i=L_{i-1}\oplus L_i$ (indices modulo $4$).
Indeed, 
the fact the Stokes gradings split the Stokes filtrations 
implies the ``subdominant'' line in each Stokes grading equals
the subdominant line in the Stokes filtration.
Then the fact the Stokes gradings are continuous across each oscillating direction $\IS\subset \partial$ implies the ``dominant'' lines in each Stokes grading are also determined.
The general case is trickier.
\section{Canonical splittings}\label{sn: canl splittings}

For the existence of such splittings, the approach here fleshes out a sketch of Malgrange 
\cite{DMRci} p.73, composing one level splittings
(this statement can be obtained in several other ways, e.g. via multisummation).
This is then upgraded to 
give the uniqueness statement in Thm. \ref{thm: mainresult} too, 
which looks to be new.
The idea of decomposing by the levels reflects the Gevrey asymptotics \cite{ramis-factn, ramis-pv}, 
and was previously used in the algebraic construction  of the 
general (multi-level) wild character varieties 
\cite{gbs, twcv}, nesting the one-level fission spaces.

\subsection{Facts about Stokes groups etc.}\label{sn: facts about stokes groups}

This section collects some facts about Stokes groups and their relation to splittings and the splitting
groups $\La_d$.
Similar statements are in \cite{BJLproper}, \cite{MR91} p.362, \cite{L-R94} (some terminology comes from the $G$ version in \cite{bafi, gbs}).

Let $(V,\Th,F)$ be a Stokes filtered local system with active exponents $I\to \partial$.
Let $V^0=\Gr(V,F)\to \partial$ be the
associated  $I$-graded local system.

Let $\ISto_d \subset \GL(V_d^0)$ be the 
Stokes group, determined by the Stokes arrows 
$\prec_d$ as in \eqref{eq: st gp}.
It is trivial unless $d\in \IA$.

Similarly define the splitting group
$\La_d \subset \GL(V_d^0)$ 
to be the unipotent group with Lie algebra determined
by the dominance (partial) order $<_d$.
This is maximal if $d\not\in \IS$, and if $d\in \IS$ it is the intersection of the two adjacent (maximal)
groups.
The set  $\Splits_d\subset \Iso(V^0_d,V_d)$ of 
splittings of $(V_d,F_d)$ is a torsor for $\La_d$.

More generally for $\al,\be\in \IR_{>0}$ 
let $\Splits_d(-\al,\be)\subset \Splits_d$
be the (possibly empty) 
set of gradings that split the 
Stokes filtrations throughout $\Sect_d(-\al,\be)$
(recalling the conventions from \S\ref{ext ints} for ramified intervals). 
If nonempty, this is a torsor for the group 
$\La_d(-\al,\be)\subset \La_d$ (defined by transporting to $d$ all the splitting groups over 
$\Sect_d(-\al,\be)$ and  then intersecting them).
Similarly for closed/half-closed intervals etc (if a closed interval ends on a Stokes direction then the intermediate filtration there is included).

Now suppose $I$ just has one level $k$.
Then $\IS = \IA+\pi/2k = \IS+\pi/k$, and this leads to a clear dictionary between the Stokes groups and the splitting groups, as follows. 
All this follows directly from the definitions of $\ISto_d$ and $\La$.

First some basic terminology:
1) a singular sector $U\subset \partial$
is an open interval bounded by consecutive singular directions, 
2) a half-period $\bd$ is a sequence of consecutive singular directions turning in a positive sense, 
consisting of the singular directions under some 
(possibly ramified) sector of the form 
$\Sect_d[0,\pi/k)$, 
3) the supersector $\wh U$ of any singular sector 
$U$ is the sector
$\wh U = \Sect_U(-\pi/2k,\pi/2k)$ (recalling 
\eqref{eq: sectU}).
Note $\wh U$ is bounded by Stokes directions and 
that the underlying singular directions 
make up a half-period. 
4) If $\bd=(d_1,\ldots,d_l)$ is a half-period, let 
$d_0=d_l-\pi/k\in \IA$ and let 
$\th(\bd):= (d_0+d_1)/2 + \pi/2k$ be the 
bisecting direction  of $\bd$ (it is not a Stokes direction).

The ``restrictions'' of $I,V,V^0$ to 
$\Sect_d(-\al,\al)$ (as in 3 of \S\ref{ext ints})
are trivial 
so it makes sense to compare the 
partial orders $\prec_d,<_d$ (and thus the 
Stokes groups and the splitting groups) defined at different points in such ramified sectors (using the transport in the sector). 

The key fact is that 
$\prec_d$ is extended by $<_e$ if 
$e\in \Sect_d(-\pi/2k,\pi/2k)$.
This implies 
$$\ISto_d\subset \La_e
\qquad\text{for all\ } e\in \Sect_d(-\pi/2k,\pi/2k).$$
Moreover any relation $i<_d j$ corresponds uniquely to a relation $i\prec_e j$ for a unique
$e\in \Sect_d(-\pi/2k,\pi/2k)$.
This implies, for any $d\in\partial$, that
$\La_d$ is directly spanned by the 
Stokes groups $\ISto_e$ for 
$e\in \Sect_d(-\pi/2k,\pi/2k)$:
$$\La_d=\left\< \, \ISto_e\st e\in \Sect_d(-\pi/2k,\pi/2k) \, \right\>^\oplus.$$
Recall (\cite{Bor91} \S14)
that a group is ``directly spanned''
by a collection of subgroups if the product map (in any fixed order) is an isomorphism (of spaces not necessarily of groups).
For example the Stokes groups in any half-period
$\bd$
directly span a full unipotent group (the unipotent radical of a parabolic)
$\ISto_\bd := 
\left\< \, \ISto_d \st {d\in \bd} \right\>^\oplus\subset \GL(V^0_{\th(\bd)}).$
It follows that $\La_{\th(\bd)} = \ISto_\bd.$ 

More generally if $\al+\be<\pi/k$ then
 $\La_d[-\al,\be]$  is
directly spanned by the Stokes groups 
in $\Sect_d(\be-\pi/2k,\pi/2k-\al)$:
\beq\label{eq: decomp of Stokes torsor}
\La_d[-\al,\be]=\bigcap_{e\in\Sect_d[-\al,\be]}\La_e = 
\left\< \,  \ISto_f \st {f\in \Sect_d(\be-\pi/2k,\pi/2k-\al)}\, \right\>^\oplus.
\eeq

Conversely if $d\in \IA$ let 
$\bd^+$ be the half-period of singular directions in
$\Sect_d[0,\pi/k)$ and let 
$\bd^-$ 
be the half-period of singular directions in 
$\Sect_d(-\pi/k,0]$ (the two 
half-periods ending on $d$).
Then
\beq\label{eq: stod from la}
\ISto_d = \ISto_{\bd^-}\cap \ISto_{\bd^+}=
\La_d(-\pi/2k,\pi/2k).
\eeq

\subsection{Preferred splittings in the one level case}\ 
If $d\not\in \IA$ then \eqref{eq: stod from la} still holds, and now
says that $\La_d(-\pi/2k,\pi/2k)=\{1\}.$
If $U$ is the singular sector containing $d$
this is the same as saying 
$\La_d(\wh U)=\{1\}$, where $\wh U$ is the supersector. 
Note that $d\pm \pi/2k$ are not Stokes directions---since they
differ by $\pi/k$ their (total) dominance orderings  
are
{\em opposite}, and this is the reason the splitting group is trivial.
This implies 
there is at most one splitting on $\wh U$.
A key point is that there is exactly one:
\begin{prop}\label{prop: one level splittings}
If $(V,\Th,F)$ has just one level $k$ and 
$U\subset \partial$ is a singular sector
then there is a unique splitting  on
the supersector 
$\wh U =\Sect_U(-\pi/2k,\pi/2k)$.
Equivalently if
$d\in \partial\setminus \IA$ 
then  
$\Splits_d(-\pi/2k,\pi/2k)$ consists of exactly one point.
\end{prop}
\pf
A direct proof is in \cite{Mal-irregENSlong} Lemme 5.1 
(see also \cite{BJL79}).
The result also follows for general reasons:
From the opposite filtrations at $d\pm \pi/2k$ it is easy to see what the grading must be.
Then one needs to check this grading splits the intermediate filtrations. 
However the Bruhat decomposition implies the intermediate
filtrations are uniquely determined by the filtrations at the ends  
(since the permutations at the Stokes directions give a length preserving decomposition of the order reversing permutation),
cf. e.g. \cite{del-lusz76} p.106 $(a_2)$.
\epf

Thus given a singular sector $U$ there is a unique splitting on $\wh U$. 
Let $\Phi_U\in \Splits(U)$
be the restriction of this splitting to $U$.
These yield a Stokes grading, and it is unique:

\begin{prop}\label{prop: one level Stokes grading}
Suppose 
$(V,\Th,F)$ is a Stokes filtered local system 
indexed by $I$, and $I$ has just one level.
Then
there is unique Stokes graded local system 
$(V,\Th,\Ga)$ that is compatible with $F$.
\end{prop}

\pf
Let $k$ denote the level, and let $V^0\to \partial$ be the associated graded.
For each singular sector $U$
let $\Phi_U\in \Splits(U)\subset \Iso_U(V^0,V)$ 
be the preferred splitting, from Prop. 
\ref{prop: one level splittings}.
To check that the Stokes conditions hold, note that if $d\in \IA$ and $U_1,U_2$ are the adjacent singular sectors then 
$\wh U_1\cap \wh U_2 = \Sect_d(-\pi/2k,\pi/2k)$.
Thus both splittings work on this intersection so they are related by an element of 
the splitting group 
$\La_d(-\pi/2k,\pi/2k)$. 
By \eqref{eq: stod from la} this group 
equals  $\ISto_d$, so the Stokes conditions hold.

Now to prove uniqueness,
let $U_0$ be a singular sector and let
$U_1,U_2,\ldots,U_{m-1}=U_{-1}$ 
be the others, in a positive sense 
(indices modulo $m$).
Let $d_i$ be the negative edge of $U_i$ so 
$U_i=\Sect(d_i,d_{i+1})$. 
Suppose $\Phi_i$ is a splitting on $U_i$ for each $i$, and they satisfy the Stokes conditions.
We will show that
$\Phi_0$ extends to a splitting on all of
$\wh U_0$, so $\Phi_0$ is the restriction of the 
unique splitting there.

Write $d=d_0,e=d_1$ so $U_0=\Sect(d,e)$ and
 $\wh U_0 = \Sect(d-\pi/2k,e+\pi/2k)$.
Let $\tau_1,\ldots, \tau_r=e+\pi/2k$ be the Stokes directions in $[e,e+\pi/2k]$, turning in 
a positive sense from $U_0$.
It follows that $\tau_{r-1} = d+\pi/2k$.
We will show inductively that
$\Phi_0$ extends to a splitting on
$(d,\tau_j)$ for $j=1,\ldots,r$.
If $j=1$ then there is nothing to do (since splittings always extend up to the next Stokes direction).
Now assume $\Phi_0$ extends to a splitting on
$(d,\tau_{j})$
for some $j<r$.
The aim is to show $\Phi_0$ extends
across $\tau_j$ to a splitting on 
$(d,\tau_{j+1})$. Observe the following:

1) there is a splitting on $(d,\tau_{j+1})$ (e.g. the unique splitting on  $\wh U_0$).

2) there is some index $d_i\in[0,\tau_j)$ such that 
$\Phi_i$ is a splitting in some (small) neighbourhood
of $\tau_j$
(by Cor. \ref{cor: malleability}, 
since $\tau_j$ will be in the closure 
of some $U_i$). 
Let $\bd^+$ be the set of 
singular directions in the sector 
$[e,\tau_{r-1})$, so $d_i\in \bd^+$
and let 
$\ISto^+ = \<\ISto_f \st {f\in \bd^+}\>^\oplus$.
By hypothesis $\Phi_i = \Phi_0 S$ for some 
$S\in \ISto^+$.

3) by definition
$\bd^+\subset (d,\tau_{r-1})=(d,d+\pi/2k)$
so that \eqref{eq: decomp of Stokes torsor}
implies 
$$\ISto^+\subset \La[d,d+\pi/2k]
= \La[d,\tau_{r-1}] =  \La[d,\tau_r)$$
where the last equality follows since $\La$ only changes 
at Stokes directions.
In particular since $\La[d,\tau_r)\subset 
\La[d,\tau_j)$ this implies
that $\Phi_0 S$ is a splitting on $(d,\tau_j)$
(since $\Phi_0$ is), 
and so $\Phi_0 S$ is in fact a splitting on all of  
$(d,\tau_{j+1})$ (as it equals the splitting $\Phi_i$ across $\tau_j$).
Thus since $S\in  \La[d,\tau_r)\subset \La(d,\tau_{j+1})$
it follows that $\Phi_0$ itself is a splitting on 
$(d,\tau_{j+1})$, completing the inductive step.

Thus, by induction, 
$\Phi_0$  is a splitting on $(d,\tau_{r})$.
Similarly going in the negative direction, 
$\Phi_0$  is a splitting on all of $\wh U_0$, and so equals the unique splitting there.
Repeating on each singular sector 
 yields the desired uniqueness statement.
\epf

\subsection{Preparation for the induction}

In order to get the proof to work cleanly a slightly more general context is needed, which will now
be made explicit so as to clarify the logic. 

{\bf Natural quotients.} \
If $k$ is a positive rational number,
let $\cI^k\subset \cI$ be the sublocal system of exponents of slope $\le k$, and let 
$\cI(k)=\cI/\cI^k$ be the quotient local system.
Thus sections of $\cI(k)$ can be represented by 
functions that can be expressed as finite sums 
of the form $q=\sum a_i x^{l_i}$ for rational numbers 
$l_i>k$ (where $x=z^{-1}$ for a local coordinate $z$). 
In turn if $I\to \partial$ is a finite subcover,
let
$$J = I(k) =  I/\cI^k.$$
This means that two local sections of $I$ 
are identified if their difference has slope $\le k$. 
Note that the maps $I\onto J$ and  $J\to \partial$ 
are both finite covering maps,
so there is a factorisation  
$$I\onto J\onto \partial$$
of the cover $I\to \partial$.
Such covers $J=I(k)\to \partial$
will be called ``natural quotients''.
If $l>k$ one can repeat and define 
$K=J(l)=J/\cI^l=I(l)$, so that  $$I\onto J\onto K\onto \partial.$$

{\bf Partitions of the fibres.} \ 
Given a nested cover $\pi:I\onto J\onto \partial$ (as above)
then $I$ can be viewed as ``graded'' 
or partitioned by $J$.
Namely each fibre $I_d$ of $I$ is partitioned into ``parts'' $\pi^{-1}(j)$, indexed by $j\in J_d$ 
(for any $d\in \partial$). 

{\bf Canonical factorisation and fission tree.} \
Thus for any natural quotient $\pi:I\to \partial$ 
(such that $\pi$ is not an isomorphism) 
there is a minimal 
rational number $k\in \IQ$ such that 
the map $I\onto I(k)$ is not an isomorphism.
Iterating, it follows  that
any finite subcover $I\to \partial$
has a canonical factorisation
\beq\label{eq: level factn}
I \onto I(k_1) \onto I(k_2) \onto \cdots \onto
I(k_{r-1}) \onto I(k_r)\mapright{\cong}\, \partial
\eeq
where each $k_i$ is minimal so that
the map $I(k_{i-1}) \onto I(k_{i})$ is not an isomorphism.
The numbers $k_1<k_2\cdots <k_r$ are the 
levels (i.e. the adjoint slopes) of $I$.
In particular 
the sequence of degrees of the covers 
$I(k_i)\to \partial$ is strictly decreasing.
Sometimes this nested sequence of covers
\eqref{eq: level factn} will be called the ``three-dimensional fission tree'' of $I$.
(It is a quite remarkable topological object canonically associated to any algebraic connection.)
If $I$ is unramified then each map $I(k)\to \partial$ 
is a trivial finite cover, and then
the three-dimensional fission tree is just the product of the circle $\partial$ with the (two-dimensional) 
{fission tree} 
attached to $I$
(the rooted tree 
described in \cite{rsode} Apx. C).
Indeed if we quotient by the action rotating the circle
a tree becomes visible: 
if a coordinate is chosen,
$I/\partial$ can be identified with a finite subset
$\{q_1,\ldots,q_n\}$
of $x\IC[x]$ (or more precisely as a multiset).
In turn $(I(k)/\partial)\subset x\IC[x]$ 
can be defined by deleting all 
monomials of each $q_i$ of degree $\le k$, 
yielding the map $I(k)\onto I(l)$ for any 
$l>k$.
Thus the quotient of 
\eqref{eq: level factn} by the circle is visibly a tree
with leaves $I/\partial$ and root 
$I(k_r)/\partial =\{*\}$:
\beq\label{eq: 2d level factn}
I/\partial \onto I(k_1)/\partial \onto \cdots  \onto I(k_r)/\partial =\{*\}.
\eeq
The picture on the title page corresponds to a simple 
twisted case (related to the irregular class 
$\<x^{5/2}\>$ at $\infty$ of the linear equation whose isomonodromic deformations yield the first Painlev\'e equation).
The reader is invited to draw (or at least imagine) the analogous picture in a case with several levels (iterated fission).

{\bf Data attached to natural quotients.} \
Natural quotients $J=I(k)$ behave just like 
finite subcovers of $\cI$. 
In particular  
the dominance ordering descends to the fibres of $J$
in the obvious way and then  Stokes directions, singular directions, levels, Stokes arrows etc of $J\to \partial$
are well-defined.
In turn the notion of 
``Stokes filtered local system indexed by $J$''
is well-defined. 
Indeed we could choose 
a component of $\cI$ 
in each equivalence class so as to
get an embedding
$J \into \cI$, and thus identify 
$J$ itself as a finite subcover
(for example by forgetting all 
the terms of order $\le k$ with respect to some coordinate).
Then all the data attached to $J$
is the same as that which arises by viewing $J$ as a 
finite subcover (independent of the choice of embedding). 
In this way the notion of a natural quotient is a mild generalisation of a finite subcover of $\cI$ (since we don't want to choose such an embedding), and we don't need to worry about the fact that $J$ is not 
canonically a finite subcover of $\cI$.

{\bf Partially graded Stokes filtered local systems.} \ 
Suppose $V\to \partial$ is a $J$-graded local system.
Then the dominance ordering determined by $J$
determines the associated filtration, 
yielding a Stokes filtered local system $(V,\cF(V))$
indexed by $J$. 
This will be called the tautological Stokes filtration.

Now suppose $\pi:I\to J$ is a natural quotient.
A ``Stokes filtered local system indexed by $I\to J$'' is a  pair $(V, F)$ where 
$V$ is a $J$-graded local system
and $(V,F)$ is a Stokes filtered local system
indexed by $I$, and they are compatible in the sense that each $I$-filtration $F_d$ refines the 
$J$-filtration $\cF(V)_d$.
In particular each local $I$-grading 
$\bigoplus_{I_d} V_d(i)$
(in the definition of $(V,F)$), 
is a refinement of the $J$-grading:
$$V_d(j) = \bigoplus_{i\in \pi^{-1}(j)} V_d(i)$$
for all $j\in J_d$.
For clarity this will sometimes be called a ``partially graded Stokes filtered local system''. 
For example if $J=\partial$ then this is just a 
Stokes filtered local system indexed by $I$.
At the other extreme, if $I=J$ it is just an $I$-graded local system (with its tautological Stokes filtration).
Some natural examples will appear in the next subsection.

Since the exponents in different graded pieces do not interact it is natural to define the 
``levels of $I\to J$'' 
to be the slopes of local sections
$q_i-q_j$ where $i,j$ are in the same part of $I$.
This means that  the levels of $I\to I(k)$ 
are the levels of $I$ that are $\le k$.
Thus in the factorisation \eqref{eq: level factn}
each map 
$I(k_{i}) \to I(k_{i+1})$ just has one level, equal to  $k_i$
(including $I\to I(k_1)$).
In turn the Stokes directions, singular directions, Stokes arrows etc of $I\to J$
are well-defined.
For example the Stokes arrows 
$I\!\to\! J$ are the subset of $I\times I$
given by the pairs $(q_i,q_j)$ where $i,j$ are in the same part of $I_d$ and $d$ is a point of maximal decay for $q_i-q_j$.
The corresponding singular directions are 
denoted $\IA(I\!\to\!J)$.

Similarly a ``Stokes graded local system 
$(V,\Ga)$ indexed by $I\!\to\! J$'', consists of 
a  $J$-graded local system $V$, plus an $I$-grading 
$\Ga_d$ of $V_d$ that refines the $J$ grading, for each $d\in \partial\setminus\IA(I\to J)$, and satisfies the Stokes conditions.

If $I\!\to\! J$ just has one level and 
$(V,F)$ is indexed by $I\to J$
then 
Prop. \ref{prop: one level Stokes grading} 
cannot be applied blindly  
to see there is a unique 
compatible Stokes grading,
since 
1) $J$ may have monodromy, and 
2) the dominance orders of $I_d$ 
are not necessarily opposite at the ends of a 
supersector (but they are in each part).
However it is clear preferred gradings exist:
given a singular sector $U$,
consider one sheet of $J$ on $U$ and the corresponding 
part of $I$---upon restriction to the 
supersector $\wh U$,
Prop. \ref{prop: one level splittings}
can then be applied to give a unique splitting.
Then repeat for each sheet of $J$.
Then the proof of Prop.
\ref{prop: one level Stokes grading} 
works verbatim to show this gives a Stokes grading and that there are no others, yielding:
\begin{prop}\label{prop: one level Stokes grading v2} 
Suppose $I\!\to\! J$ just has one level and 
$(V,F)$ is a Stokes filtered local system 
indexed by $I\to J$.
Then
there is unique Stokes graded local system 
$(V,\Ga)$ indexed by $I\to J$,
that is compatible with $F$.
\end{prop}

{\bf Stokes groups by level.} \ 
Suppose $I\to J=I(k_i)$.
If $V$ is $I$-graded 
and $d\in \partial$
then the  three Stokes groups:
$$\ISto_d(I),\   \ISto_d(J), \ \ISto_d(I\!\to\!J)\subset 
\GL(V_d)$$
are well defined.
\begin{lem}\label{lem: vertical directspanning}
The groups $\ISto_d(J), \ISto_d(I\!\to\! J)$
are subgroups of $\ISto_d(I)$ 
(with $\ISto_d(J)$ being a normal subgroup)
and they directly span it in any order.
In other words there is a semidirect product decomposition:
\beq\label{eq: 2 step level decomp}
\ISto_d(I) = \ISto_d(I\to J)\sdp \ISto_d(J).
\eeq
\end{lem}
Indeed there is a parabolic subgroup 
$P\subset \GL(V_d)$ with Levi decomposition
$P=H\cdot U$ such that 
$\ISto_d(I)\subset P$ and 
$\ISto_d(I\to J) = \ISto_d(I)\cap H$ and 
$\ISto_d(J) = \ISto_d(I)\cap U$.
Here $H$ is the automorphism group of the $J$-grading 
of $V_d$ (not the $I$-grading).

Iterating as much as possible (as in \eqref{eq: level factn}) this gives the level decomposition of the Stokes groups, i.e. the direct spanning deomposition:
\beq\label{eq: full decomp}
\ISto_d(I) = \<\ \ISto_d^k(I)\st \text{ $k$ is a level of $I$} \ \>^\oplus
\eeq
where $\ISto_d^k(I)$ is the level $k$ Stokes group, i.e. the subgroup corresponding to the level $k$ Stokes arrows, i.e. the arrows  $i\prec_d j$ such that $q_i-q_j$ has slope $k$.
The decomposition
\eqref{eq: 2 step level decomp} corresponds to taking the arrows of level $\le k_i$ (to give $\ISto_d(I\to J)$) and the arrows of level $>k_i$ 
(to give $\ISto_d(J)$).

{\bf Partial associated gradeds.}\ 
Suppose $(V,F)$ is a Stokes filtered local system indexed by $I\to K$ with $K=I(l)$.  
Choose $k<l$ and let $J=I(k)$ so that 
$I\to J\to K\to \partial$.
Then two new Stokes filtered local systems
indexed by $J\to K$ and by $I\to J$ respectively, can be defined, as follows.

First define a filtration $F^{(k)}$ 
indexed by $J$ on $V$:
$$F^{(k)}_d(j) = \sum_{i\in j} F_d(i)\qquad 
(\text{where } j\in J_d
\text{ so } j\subset I_d).$$
\begin{lem}
$(V,F^{(k)})$ is a Stokes filtered local system indexed by $J\to K$.
\end{lem}
\pf
Given a local splitting (grading) 
of $(V,F)$ indexed by $I$,
we can collapse the grading to be indexed by $J$.
Then $F^{(k)}$ is the associated filtration 
indexed by $J$.
\epf

Now let $V^k = \Gr(V,F^{(k)})$ be the corresponding 
associated graded local system.
It is a $J$-graded local system.
Let $F^k$ be the filtration 
on the fibres of $V^k$ induced by $F$.
\begin{lem}
$(V^k,F^k)$ is a %
Stokes filtered local system indexed by $I\to J$, and its associated graded is canonically 
isomorphic to $\Gr(V,F)$ (as $I$-graded local systems).
\end{lem}
\pf
We will drop the point $d\in\partial$ from the notation.
By definition 
$V^k(j) = F^{(k)}(j)/F^{(k)}(<j) = 
\sum_{i\in j} F(i)/\sum_{i< j} F(i)$
where $i\in I$ and  $i<j$ 
means there is some $j'<j$ with $i\in j$.
Write $\pi_j: F^{(k)}(j)\onto V^k(j)$ 
for the natural projection.
Then for all $i\in I, j\in J$:
$$F^k(i) = \bigoplus_{j\in J} F^k(i)\cap V^k(j) \quad\text{where}\quad F^k(i)\cap V^k(j) = \pi_j( F(i)\cap F^{(k)}(j)).
$$
As above, a local splitting (grading) 
of $(V,F)$ indexed by $I$,
will induce a local splitting 
of $(V,F^{(k)})$, and thus a local 
isomorphism with its associated graded
$(V^k,\cF(V^k))$
(with its tautological filtration indexed by $J$).
This gives a local $I$-grading of $V^k$ and $F^k$ is the associated filtration, thereby showing it satisfies the Stokes conditions (when we start with a local grading across a Stokes direction).
The isomorphism of the associated gradeds is left as an exercise.
\epf

{\bf Completing the proof.}\ 
The main result, 
Thm \ref{thm: v2 of main thm}, is a special case of the more general statement:
Any Stokes filtered local system $(V,\Th,F)$ 
indexed by $I\to K$
admits a unique Stokes grading $\Ga$, such that 
$F$ is the Stokes filtration associated to $\Ga$.

In turn this can be proved by induction on the number of levels of $I\to K$.
The one level case is Prop. \ref{prop: one level Stokes grading v2}.
If $(V,F)$ 
has more than one level choose $J=I(k)$ so that 
$I\to J\to K$ and both of $I\to J$ and $J\to K$ 
have fewer levels.
Then $(V,F^{(k)})$ and $(V^k,F^k)$ both have unique 
Stokes splittings by induction.
Composing these gives a splitting of $(V,F)$ on each singular sector. The level decomposition 
\eqref{eq: 2 step level decomp},\eqref{eq: full decomp}
of the Stokes groups implies these are a Stokes grading.
Again by the decomposition of the Stokes groups, any other Stokes splitting will give a splitting of the two pieces, and thus equal the previous one (by the uniqueness of the lower level splittings).
This completes the proof.

It is visually helpful to draw the punctured disk model of the Stokes local system, decomposing each tangential puncture radially into several tangential punctures 
(one for each level, supported by the corresponding components of the Stokes groups), as in 
\cite{gbs}\S7.2. Then the direct spanning decompositions \eqref{eq: 2 step level decomp},\eqref{eq: full decomp}
amount to factorising loops around these decomposed tangential punctures (this amounts to nesting one level fission spaces as in \cite{gbs, twcv}, and reflects the Gevrey filtration \cite{ramis-factn, ramis-pv} as already mentioned).

\subsection{Fully faithfullness}\label{sn: ff}
Here we will prove fully faithfulness, 
which hinges on the following 
proposition.
Recall that $\<0\>\subset \cI$ is the tame circle. 
\begin{prop}\label{prop: fflem}
Let $(V,F)$ be a Stokes filtered local system
indexed by $I\to \partial$
and let $\Ga$ be the 
unique compatible Stokes grading.
Suppose $v$ is a section of $V$ on $\partial$,
and $v(d)\in F_d\<0\>$ for all 
$d\in \partial\setminus \IS$.
Then $v(d)\in \Ga_d\<0\>$ for all 
$d\in \partial\setminus \IA$.
Consequently (by Lemma \ref{lem: monod invt vectors}) 
$v$  extends uniquely to a section of $\Gr(V,F)_{\<0\>}$
(the tame piece of the associated graded local system).
\end{prop}

The following corollary shows how the topological picture encodes the fact  that if a solution has exponential decay in some direction 
then it will have exponential growth somewhere else 
(moderate global solutions cannot decay exponentially
anywhere).

\begin{cor}
Suppose $v,(V,F)$ satisfy the hypotheses of Prop. 
\ref{prop: fflem}. If there exists 
$d\in \partial\setminus \IS$ and $i\in I_d$  such that
 $i<_d \<0\>$ and $v(d)\in F_d(i)$,  then $v=0$.
\end{cor}
\pfms(of corollary).
One of the gradings $\Ga$ will split $F$ across
$d$.   
\epfms

\begin{cor}\label{cor: equal spaces on bdy}
The following two subspaces of 
$\HH^0(V,\partial)$ are equal:

$\HH^0((V,\Ga),\partial) := \{ v\in \HH^0(V,\partial)\st 
v(d)\in \Ga_d\<0\>
\text{ for all } d\in \partial\setminus \IA\},$

$\HH^0((V,F),\partial) := \{ v\in \HH^0(V,\partial)\st v(d)\in F_d\<0\>
\text{ for all } d\in \partial\setminus \IS\}$

\noindent
and moreover this common vector space
embeds naturally in the space
$\HH^0(\Gr(V,F)_{\<0\>},\partial)$
of sections of the tame component of the associated graded local system. 
\end{cor}
\pf
It is clear that the first space is contained in the second.
The reverse inclusion follows from 
Prop. \ref{prop: fflem}, as does the inclusion
in $\HH^0(\Gr(V,F)_{\<0\>},\partial)$.
\epf

\begin{cor}
The functor $\varphi$ taking a Stokes graded 
local system to the associated Stokes filtered local system is fully faithful.
\end{cor}
\pf
Using the internal hom, 
this comes down to showing that there is an equality 
$\HH^0((V,\Ga),\wh \Si) = \HH^0((V,F),\wh \Si)$
of subspaces of $\HH^0(V,\wh \Si)$, 
for any Stokes graded local system
$(V,\Ga)$, where $(V,F)=\varphi(V,\Ga)$.
This follows from Prop. \ref{prop: fflem} as in Cor. \ref{cor: equal spaces on bdy}.
\epf

\pfms(of Prop. \ref{prop: fflem}).
First replace $I$ by $I\cup \<0\>$ (to avoid a separate argument ruling out the case with $\<0\>$ absent---it will follow that $v$ is zero then).
We will prove (by induction on the number of levels) the more general statement with $I\to \partial$ replaced by 
$I\to J$.

First suppose $I\to J$ just has one level (and $\<0\>$ is present in both $I$ and $J$).

Its enough to prove that if $v(d)\in F_d(i)$ for some 
$d\in \partial\setminus \IS$ and $i\in I_d$ with 
$i<_d \<0\>$, then $v=0$.
(Indeed if this holds then 
$v\in F_d\<0\>$ will imply $v\in \Ga_d\<0\>$.)

Write 
$\partial \setminus \IS = 
U_0\cup U_1 \cup \cdots\cup U_{m-1}$ 
and suppose $d\in U_0$.
Choose $i<_d \<0\>$ minimal in $I_d$ 
so that $v\in F_d(i)$.
Let $c=c(i)\in \IN$ be the number of Stokes directions
one needs to cross (in a positive sense) before the dominance ordering between $i$ and $\<0\>$ changes (possibly going around the circle several times).
If we pass to $U_1$ and $i$ does not cross $\<0\>$,
possibly another index $j$ may replace $i$ (as the 
minimal index so that $v\in F_d(j)$),
but in any case the number 
$c$ will decrease by at least $1$:

\begin{lem}
If $j<_d i<_d \<0\>$ in $U_0$ and   
$i<_d j<_d \<0\>$ in $U_1$ then $c(j)\le c(i)$.
\end{lem}
\pf
Since there is only one level $k$ the dominance order of each pair in $i,j,\<0\>$ changes every $\pi/k$. 
Thus if $j$ crosses above $i$ then $j$ will cross $0$ before $i$ does.
\epf

Thus after a finite number of steps there is a Stokes direction $\tau$ where the minimal index $i$ 
will change dominance with $\<0\>$.
Choosing a local splitting of the Stokes filtrations across $\tau$, 
we see that if $v$ is nonzero then $v\not\in F\<0\>$
on the other side of $\tau$.
Thus $v$ is zero.

The general statement can now be deduced.
If $(V,F)$ is indexed by $I\to K$ with $>1$ levels, 
then as before there is a factorisation 
$I\to J \to K$, and $(V,F^{(k)})$ indexed by 
$J\to K$, and $(V^k,F^k)$ indexed by $I\to J$, both with fewer levels.
Given a section $v$ of $V$ with 
$v(d)\in F_d\<0\>$ 
then $v(d)\in F^{(k)}_d\<0\>$,
so we can apply the inductive hypothesis to
$(V,F^{(k)})$.
Thus $v$ takes values in the tame piece of the 
Stokes grading 
$\Ga^{(k)}$ of $V$ (that splits $F^{(k)}$), and extends uniquely to a section of the tame component of the 
associated graded, i.e. $V^k\<0\>$.
Thus we get a section of $V^k$, and by the original 
hypothesis it lives in the piece $F_d^k\<0\>$ of
the induced filtration.
Thus we can apply the inductive hypothesis to
$(V^k,F^{k})$ and see that $v$ takes values
in the tame piece of its Stokes grading $\Ga^{k}$ 
and extends uniquely to a 
section of the tame component of the
associated graded, 
$\Gr(V^k,F^k)=\Gr(V,F)$.
This is the desired statement 
(since $\Ga^k$ gives $\Ga$ once we use view 
the splitting $\Ga^{(k)}$ as giving a local isomorphism
$V^k\cong V$).
\epfms

\section{Wild character varieties and moduli problems}\label{sn: wcvs}

This section will review the main implications for the 
wild character varieties (considered in the generic case in \cite{birkhoff-1913, JMU81, smid} and in general
in \cite{gbs, twcv}, using the local theory of 
\cite{L-R94}). 
Fix a smooth complex projective curve $\Si$ and some
marked points $\ba\subset \Si$.
Write $\Si^\circ = \Si\setminus \ba$.

\subsection{Tame character varieties}
This section will quickly run through the
theory of tame character varieties, as a model before discussing the wild case.
Let $\LocSys$ be the category of local systems of finite dimensional complex vector spaces
(on the topological surface underlying $\Si^\circ$).
Each such local system has an invariant, its rank.
Let $\LocSys(n)$ be the groupoid of local systems of rank $n$ (so that isomorphisms are the only maps considered).

The set of isomorphism classes in $\LocSys(n)$
appears as the set of orbits of a complex reductive group on an affine variety, as follows.  
Choose a basepoint $b\in \Si^\circ$ then 
$\LocSys(n)$ is equivalent to the category of 
rank $n$ modules for the group $\pi_1(\Si^\circ,b)$ 
(taking a local system to its monodromy representation).
Recall that a framed local system is a local system $V$ equipped with a framing at $b$, i.e. 
a basis  $\phi:\IC^n\to V_b$ of the fibre at $b$.
The {\em representation variety}
$$\cR_n=\Hom(\pi_1(\Si^\circ,b), \GL_n(\IC))$$
is the set of isomorphism classes of framed 
rank $n$ local systems.
Choosing a presentation of $\pi_1(\Si^\circ,b)$
makes it clear that $\cR_n$
is an affine variety (replace the generators by elements of 
$\GL_n(\IC)$ in the relation).
Changing the framing corresponds to the 
conjugation 
action of $G=\GL_n(\IC)$ on $\cR_n$, and it follows that
the 
set of isomorphism classes in $\LocSys(n)$
is in bijection with the set of $G$ 
orbits in the affine variety $\cR_n$.

The {\em character stack} $\gM_n$ is the stack theoretic quotient of $\cR_n$ by $G$, whereas 
the {\em character variety} $\cM_n$ 
is the affine geometric invariant theory 
 quotient of $\cR_n$ by $G$,
i.e. $\cM_n$ is the variety associated to 
the ring
$$\IC[\cR_n]^G$$
of $G$-invariant functions on the affine variety
$\cR_n$. This ring is finitely generated since $G$ is reductive.

It is well-known that $\cM_n$ has an algebraic Poisson structure with symplectic leaves given by fixing the 
isomorphism class of the local system in a small punctured disk around each puncture 
(this is the same as fixing the
conjugacy class of monodromy around each puncture).

\subsection{Wild character varieties}
The definition of the wild character varieties
now follows a similar pattern. 
Let $\SLocSys$ be the category of Stokes local systems associated to $\Si,\ba$.
To simplify the presentation suppose
$\ba$ consists of just one point.
The general case is in \cite{gbs, twcv}.

Each such Stokes local system has an invariant, its irregular class.
Let $\SLocSys(\Th)$ be the groupoid of Stokes 
local systems of class $\Th$ (so that isomorphisms are the only maps considered).
In particular this fixes its rank, $n=\rk(\Th)$.

The set of isomorphism classes in $\SLocSys(\Th)$
appears as the set of orbits of a complex reductive group on an affine variety, as follows. (This will then yield the wild character variety $\cM_\Th$ of the irregular curve $\bSi=(\Si,\ba,\Th)$.)

Let $\wt \Si\subset \wh \Si$
be the auxiliary curve
determined by the irregular class $\Th$ (removing a tangential puncture near each singular direction).
Choose a basepoint $b\in \partial$ 
then
$\SLocSys(\Th)$ embeds in the category of 
rank $n$ modules for the group 
$\Pi:=\pi_1(\wt \Si,b)$
(the {\em wild surface group}). 
A framed Stokes local system is a Stokes local system 
$\IV$ equipped with a framing at $b$, i.e. 
an isomorphism $\phi:\IF\to \IV_b$ 
of graded vector spaces, where
$\IF=\IC^\Th:= \bigoplus_{i\in \cI_b} \IC^{\Th(i)}$
is the standard fibre---a graded vector space of dimension $\Th$,
so that $\IF(i)=\IC^{\Th(i)}$.

The {\em naive wild representation variety} is the space
$\Hom(\Pi, G)$   where $G=\GL(\IF)\cong \GL_n(\IC)$.
Each framed Stokes local system canonically determines a point of
$\Hom(\Pi, G)$
and it is easy to characterise the subvariety
of Stokes representations (the 
{wild representation variety) 
$$\cR_\Th = \Hom_\IS(\Pi, G)\subset \Hom(\Pi, G)$$
where the monodromy of Stokes local systems lives 
(reflecting the facts that the monodromy around 
$\partial$ should be the 
monodromy of a graded local system, and the 
Stokes conditions, that monodromy around each tangential puncture should land in the corresponding Stokes group). 
This will be spelt out below in \S\ref{ssn: streps subvr}.

In this way $\cR_\Th=\Hom_\IS(\Pi, G)$
is the set of isomorphisms classes of framed 
Stokes local systems of class $\Th$.
Choosing a presentation of $\Pi$
makes it clear that $\cR_\Th$
is an affine variety 
(since the Stokes groups are affine too).
Changing the framing corresponds to the 
conjugation 
action of $H=\GrAut(\IF)\subset G$ on $\cR_\Th$, and
it follows that the 
set of isomorphism classes in $\SLocSys(\Th)$
is in bijection with the set of $H$ 
orbits in the affine variety $\cR_\Th$.
Thus this key fact still persists in the wild case.

Similarly a framing of a Stokes filtered local system
$(V,F,\Th)$ is a graded isomorphism 
$\phi:\IF\to \Gr(V,F)_b$, and a 
framing of a Stokes graded local system
$(V,\Ga,\Th)$ is a graded isomorphism 
$\phi:\IF\to (V,\Ga_b)$, taking the median grading if $b$ is a singular direction.
The main result of this paper then implies the 
following statement:

\begin{cor} 
The wild representation variety
$\cR_\Th=\Hom_\IS(\Pi,G)$ 
parameterises the set of isomorphism classes of:

$\bullet$ Framed Stokes filtered local systems of class $\Th$, and also of

$\bullet$ Framed Stokes graded local systems of class $\Th$, and also of

$\bullet$ Framed Stokes local systems of class $\Th$.
\end{cor}

Note this corollary
also follows from the local result in \cite{L-R94}
after some identifications.

In turn the {\em wild character stack} $\gM_\Th$ is 
the stack theoretic quotient of $\cR_\Th$ by $H$, whereas 
the {\em wild character variety} $\cM_\Th$ 
is the affine geometric invariant theory  
quotient of $\cR_\Th$ by $H$,
i.e. $\cM_\Th$ is the variety associated to 
the ring
$$\IC[\cR_\Th]^H$$
of $H$-invariant functions on the affine variety
$\cR_\Th$. 
This ring is finitely generated since $H$ is reductive.
Stability for the action of $H$ has been analysed in 
\cite{gbs}.

Results of \cite{thesis, smid, saqh, fission, gbs,twcv} show that the wild character variety $\cM_\Th$ has an algebraic Poisson structure with symplectic leaves given by fixing the 
isomorphism class of the graded local system in the halo
(this is the same as fixing a twisted
conjugacy class for  the group $H$, cf. \cite{twcv}).
This generalises the tame case, where there are no tangential punctures, and the grading is trivial 
(everything is graded by the tame circle $\<0 \>$).

\subsection{}\label{ssn: streps subvr} 
Here is how to define the subvariety 
$\Hom_\IS(\Pi,G)\subset \Hom(\Pi,G)$.
Let $\wh\rho_\partial:I_b\to I_b$ be the monodromy of
the active exponents $I\to \partial$ 
(in a positive sense),
and let $\rho:\Pi\to \GL(\IF)$ be the monodromy of a Stokes local system $\IV$ of class $\Th$.

1) Due to the $I$-grading of $\IV$ on $\partial$, 
the  monodromy satisfies
\beq\label{eq: twist rel}
\rho_\partial(\IF(i)) = \IF(\wh \rho_\partial(i))
\eeq
for all $i\in I_b$.

If $d\in \IA$ let
$\al$ be any path in  
$\partial$ from $b$ to $d$, and let 
$\ga_d$ be the simple loop based at $d$
going out to $e(d)$ around it in a positive sense and then back to $d$.
Then let $\eta_d = \al^{-1} \circ \ga_d\circ \al$ 
be the corresponding loop based at $b$.
Let $\Sto_d\subset \GL(\IF)$ be the group obtained
by transporting $\ISto_d\subset \GL(\IV_d)$ along $\al$ to $b$ and then to $\GL(\IF)$ via the framing.  

2) The Stokes conditions then say that 
\beq
\rho(\eta_d)\in \Sto_d
\eeq
for all such $d$ and $\al$.

The subvariety $\Hom_\IS(\Pi,G)\subset \Hom(\Pi,G)$
is cut out by these two conditions.
By choosing generators of $\Pi$ this is easily made completely explicit, and yields the Birkhoff type presentations of the wild character variety, as the quotient by $H$ of the fibre at $1$ of a map
of the form
$$G^{2g} \times H(\partial)\times \Sto  \to G; \qquad
(\bA,\bB,h,\bS)\mapsto \left(\prod_1^g [A_i,B_i]\right) hS_r\cdots S_2S_1$$
where $[A,B]=ABA^{-1}B^{-1}$, 
$\Sto\subset G^r$ is a product of Stokes groups
and 
$$H(\partial)= \{h\in G\st 
h(\IF(i))=\IF(\wh \rho_\partial(i))
\text{ for all } i\in I_b \}$$
is the twist of $H$ consisting of elements satisfying the relation \eqref{eq: twist rel}. 
See \cite{gbs} equation (37) (and \cite{twcv})
for the multipoint 
case---the generic case in genus zero is in 
\cite{JMU81} (2.46),
closely related to that of Birkhoff \cite{birkhoff-1913} \S15.  
These presentations motivated the TQFT approach
to meromorphic connections
\cite{smid, saqh, fission, gbs, twcv} 
where such quotients are shown to be multiplicative symplectic quotients, and thus have natural symplectic/Poisson structures.

\begin{rmk}
The full story involves adding 
tame (Levelt--Simpson) filtrations as well, and the resulting wild character varieties will not in general be affine 
(the case here corresponds to having trivial Betti weights, denoted 
$\phi$ in \cite{hit70}, $\ga$ in \cite{wnabh}, $\be$ in \cite{Sim-hboncc}).
\end{rmk}

\begin{rmk}
In general \cite{gbs, twcv}, in order to fit well with the group-valued moment map approach, the group $\Pi$ is usually replaced by the fundamental 
groupoid $\Pi_1(\wt \Si,\be)$
where $\be\subset \wt \Si$ consists of one point in each component circle of $\partial$.
Beware that a different groupoid was used 
earlier (in \cite{smid} p.160) to encode Stokes data.
\end{rmk}

\begin{rmk}
The Stokes decompositions and 
wild monodromy/Stokes local system
don't seem to have been used in the linear setting  beyond the curve   case, 
so these approaches may well simplify and render more explicit existing approaches (cf. \cite{AgKash16}).
New examples of such wild character varieties seem lacking
(see the problem at the end of \S1.6 in \cite{hit70}).
\end{rmk}

\subsection{Wild nonabelian periods/wild Wilson loops}\label{sn: periods}

Functions on $\cR_\Th$ invariant under the action of 
$H$ may be constructed as follows.
In the tame case, one would just take 
the trace of the monodromy 
around loops in the surface.
In general one can take the wild monodromy, i.e.
the monodromy of the Stokes local system along wild loops, i.e. 
loops in $\wt \Si$. 
Moreover if the wild loop is based 
in $\partial$ then 
the fibre is graded 
so there are more invariants than just the trace
(one only needs to quotient by the graded automorphisms).

Let $K=K(I_b)$ be the complete quiver with nodes $I_b$.
This has a loop at each node and a directed edge in
each direction between each pair of distinct nodes.
Thus if $V$ is any $I_b$-graded vector space then
$\End(V)$ is the same as the space $\Rep(K,V)$ 
of quiver representations of $K$ on $V$:
$$\End(V) = \Rep(K,V).$$
Any cycle (i.e. a loop) $C$ in $K$
determines an $H$-invariant function
$$\phi_C:\Rep(K,V) \to \IC$$
by taking the trace of the composition of the 
maps along the edges in the cycle $C$.
Here   $H=\GrAut(V)$.
Now given $\rho\in \cR_\Th=\Hom_\IS(\Pi,\GL(\IF))$
and a loop $\ga\in \Pi$, i.e. a loop in $\wt \Si$ based at $b$ then 
$\rho(\ga)\in \GL(\IF)\subset \End(\IF)=\Rep(K,\IF)$.
Thus there is an $H$-invariant function
\beq
\phi_{C,\ga}:\cR_\Th\to \IC;\quad \rho\mapsto 
\phi_{C,\ga}(\rho):=\phi_C(\rho(\ga))
\eeq
for each choice of cycle $C$ in $K$ 
and wild loop $\ga\in \Pi$.
Of course given a Stokes local system
$\IV$ with irregular class $\Th$ 
there is no need to discuss framings, and one can just work with the graded vector space 
$\IV_b$; A Stokes representation 
$\wh \rho\in \Hom_\IS(\Pi,\GL(\IV_b))$ 
is intrinsically defined and this is enough to 
construct numbers $\phi_{C,\ga}(\wh \rho)$
that depend only on the isomorphism class of $\IV$.
Of course these numbers are invariants of the corresponding connection $\nabla$ too: for each choice of loop $\ga$ and 
cycle $C$ the complex number 
$\phi_{C,\ga}(\nabla):=\phi_{C,\ga}(\wh \rho)$
is well defined. They are the ``wild nonabelian periods/wild Wilson loops'' of $\nabla$.

\appendix

\section{Analytic black boxes}\label{sn: anbbx}

For completeness this section gives a quick review of 
the analytic results
used to define the topological data (or to show they have the desired properties).
The analytic results will be presented/used as black boxes, and aren't used elsewhere in this article.
The different approaches to Stokes data
arise from the use of different analytic tools. 
One particular aim 
is to show how the abstract/intrinsic language used here 
encapsulates what most authors actually do when 
defining monodromy/Stokes data.

\subsection{Cauchy and local systems}

Let $(E,\nabla)\to \Si^\circ$ be an algebraic connection on an algebraic vector bundle $E$ on $\Si^\circ$. 
Thus $\nabla$ is an operator $\nabla:E\to E\otimes\Om^1$
which is a $\IC$-linear map of sheaves, satisfying the 
Leibniz rule $\nabla(fv)=(df)v+ f\nabla(v)$ for any local section $v$ and function $f$. 
For example
take $E=\IC^n\times \Si^\circ$ 
to be trivial and $\Si^\circ\subset \IC$ 
to be the complement of a finite number of points in the complex plane, and $\nabla = d-Bdz$ for any algebraic map $B:\Si^\circ\to \gl_n(\IC)$.
The basic existence/uniqueness theorem for systems of linear ODEs implies the following.
\begin{thm}
If $\De\subset \Si^\circ$ is a disk then the space 
$$V(\De) = \{v\in \HH^0(E^{an},\De)\st \nabla(v)=0\}$$
of analytic solutions of $(E,\nabla)$
on $\De$ is a finite dimensional complex vector space, and the map
$$V(\De) \to E_b;\qquad v\mapsto v(b)$$
to the fibre of $E$ at $b$, is a linear isomorphism for any $b\in \De$.
\end{thm}
\pf
Choose a local trivialisation of $E^{an}$ (the analytic 
vector bundle determined by $E$) on $\De$ and a coordinate $z$ on $\De$, so $V(\De)$ becomes identified with the holomorphic maps $v:\De\to \IC^n$ such that 
$dv/dz = Bv$ for some matrix of holomorphic functions $B$ on $\De$ (so that the connection is $d-Bdz$).
Now use the existence/uniqueness theorem for this system of linear ODEs.
\epf

This analytic fact defines the local system 
$V\to \Si^\circ$ of solutions of $(E,\nabla)$.

\subsection{Stokes filtrations and the local asymptotic existence theorem}\label{ssn: lae}

Given the local system $V\to \Si^\circ$
of solutions of $(E,\nabla)$, the $n$-dimensional 
vector space $V_d$ 
is well-defined for any $d\in \partial$ 
(see \S\ref{ssn: robups}).

The Stokes filtration in $V_d$ is defined by looking for the subspace of recessive solutions, i.e. those with maximal exponential decay (or least growth), in some algebraic trivialisation of $E$ across the pole.
Then quotient $V_d$ by the subspace of recessive solutions and iterate to get the filtration.
The fact that this process works, and the result is a Stokes filtration, follows from the local asymptotic
existence theorem (in turn this uses the formal classification of meromorphic connections).
In general, given $d\in \partial$ and $q\in \cI_d$
the corresponding piece 
$F_d(q)\subset V_d$
of the Stokes filtration
is 
made up of the solutions $v\in V_d$ such that
$v/\exp(q)$ has at most moderate growth at $0$ 
in some open sector containing $d$.

Choose a small disk $\De\subset \Si$ containing a marked point $a$. 
Choose a coordinate $z$ on $\De$ vanishing at $a$
and a local trivialisation of $E$,  so that 
$$\nabla = d - A,\qquad A = \sum_{-N}^\infty A_i z^i dz$$
where $A_i\in \gl_n(\IC)$ and the series is convergent.

For simplicity suppose $\nabla$ is unramified
(the general case follows by descent).
Then the formal classification implies that there 
is $\wh F\in \GL_n(\IC\flp z \frp)$ such that
$$A = \wh F[A^0] = \wh FA^0 \wh F^{-1} + 
(d\wh F) \wh F^{-1}\quad \text{where } 
A^0 = dQ + \La dz/z,\ Q=\diag(q_1,\ldots,q_n).$$ 
Here $Q$ is an irregular type, with 
$q_i\in z^{-1}\IC[z^{-1}]$, and $\La\in \gl_n(\IC)$ is 
a constant matrix that commutes with $Q$.

Thus  $\wh F$ is a formal isomorphism taking 
$\nabla^0$ to $\nabla$, where  
$\nabla^0 = d-A^0$ (often called the formal normal form).
The irregular class is given by the $\<q_i\>$ with their 
multiplicities (each is a trivial cover of the circle $\partial$).
The general result is that one can always pass to a cyclic cover $(t^r=z)$  and then get such a formal normal form  
upstairs.

The local asymptotic existence theorem (cf. \cite{Was76} Thm. 19.1)
says that any direction $d\in \partial$
has a (small) open neighbourhood $U\subset \partial$
on which there exists an analytic isomorphism
$F$ taking $\nabla^0$ to $\nabla$, 
that is asymptotic at zero to $\wh F$ in $U$.

The topological interpretation is as follows:

1) $\nabla^0=d-A^0$ is a graded connection, it breaks up into a direct sum of connections indexed by the set 
$z^{-1}\IC[z^{-1}]$ of unramified irregular classes.
Thus its solutions form a graded local system $V^0$ (on a germ of a punctured disc, or equivalently on $\partial$).

2) The local asymptotic existence theorem
gives local isomorphisms between 
$V^0$ and $V$ (on $\partial$).
Since $V^0$ is graded, this gives local gradings of $V$.
These gradings are not intrinsic, but the associated filtrations are completely intrinsic, and moreover
the condition for the existence of such local gradings splitting the filtrations gives a way to axiomatise the 
filtrations. This is Deligne's idea \cite{deligne78in} yielding the Stokes filtrations and the 
axioms for Stokes filtered local systems.

3) The associated graded 
local system $\Gr(V)\to \partial$ of the Stokes filtration
is intrinsically defined.
The choice of $\nabla^0$ determines a
graded local system $V^0$ (its solutions), and the 
choice of a formal isomorphism $\wh F$
then uniquely determines an isomorphism
$V^0\to \Gr(V)$ of graded local systems
(induced by any such local analytic isomorphism $F$).
This gives a bijection between the set of such $\wh F$
and such graded isomorphisms (the categories of formal connections and graded local systems are equivalent \cite{deligne78in}).

\begin{rmk}
If instead one chooses an open cover of $\partial$
and such a local isomorphism $F$ 
on each open set, and examines how they differ on two-fold overlaps, 
one gets to the Malgrange--Sibuya cohomological 
approach \cite{Mal79, sibuya77}.   
This amounts to taking the cohomology class classifying  the sheaf of torsors determined by the sets $\Splits_d$ of splittings. 
\end{rmk}

\subsection{Summation and preferred bases}
\label{ssn: summation and preferred bases}

The general Stokes approach comes from an apparently  stronger analytic existence theorem \cite{BBRS91}, involving multisummation in general (which generalises $k$-summation, and in turn Borel summation).

In the set-up above with $\nabla=\wh F [\nabla^0]$
this says that there are a finite number of singular directions $\IA\subset \partial$
and a preferred 
analytic isomorphism $F_U\in \Iso_U(\nabla^0,\nabla)$  
canonically determined by $\wh F$ on each singular sector $U\subset \partial$.

In particular, given the choice of a solution $w$ of $\nabla^0$ on $U$ then $\wh Fw$  is a formal solution of $\nabla$, and the theorem implies this determines a preferred solution $F_Uw$ of $\nabla$.
In other words formal solutions determine 
preferred analytic solutions.
This is the form of the result discovered by Stokes (in the examples related to the Airy and Bessel equations).
Stokes used optimal truncation rather than Borel summation, but one can check the preferred solutions 
he specified (for each formal solution) are the same
as those  one gets from Borel summation.

1) On any singular sector $U$ 
the element $F_U\in \Iso(\nabla^0,\nabla)$ 
yields a grading of the local system $V$ of solutions of
$\nabla$ (since $\nabla^0$ is graded).
The key remark to make then 
is that this is the Stokes grading, and it depends only on $\nabla$ (and not the choice of $\wh F, \nabla^0$).

2) $V,V^0$ and the gluing maps $F_U$ make up a Stokes local system with gluing maps. 

This Stokes local system is isomorphic to the canonical Stokes local system (determined by the Stokes gradings in 1) via the converse part of \S\ref{sn: sls and sgls}), once $V^0$ is identifed with the canonical graded local system via the choice of $\wh F$
(via 3) of \S\ref{ssn: lae} and Rmk \ref{rmk: same gls})

In particular the corresponding Stokes representations will be isomorphic.
This is good to know since these Stokes local systems are the ones used in practice (via choices of normal forms at each pole).

Note that the equivalence between
Stokes gradings and Stokes filtrations, implies that the analytic results of
\S\ref{ssn: lae} 
and  \S\ref{ssn: summation and preferred bases}  
are in fact algebraically equivalent.

\subsection{Intuitive way to  
understand Borel summation etc}\label{sn: borel sums}

Consider the following quite familiar statement:
\beq\begin{matrix}\label{eq: stokesinterp}
\text{\it 
Sometimes a power series determines a holomorphic}\\ 
\text{\it function outside of its domain of convergence.}
\end{matrix}
\eeq
Indeed one can just consider the power 
series at zero of the function $1/(1-x)$.
This has a pole at $x=1$ so the series has radius of convergence $1$, yet clearly the function defined by the power series can be analytically continued outside the unit disk (avoiding $1$).
Similarly the power 
series at zero for $1/\sqrt{1-x}$
defines a holomorphic function in the unit disk,
however now the branch of the 
function obtained outside the unit disk depends on which side of the point $x=1$ one takes: for example at $x=2$ one may get either sign in $\pm i$ depending on the path taken. One can readily cook up more examples with more {\em singular directions} on the unit circle, and in
turn examples with arbitrarily small radius of convergence.

The key point (to intuitively 
understand Borel summation etc) 
is that the statement 
\eqref{eq: stokesinterp} may hold
{\em  even if the series has radius of convergence zero}.
Indeed this is precisely what Borel summation does: away from the singular directions
a holomorphic function is determined, and something like a different branch of the same function appears if a different direction is used,  on the other side of  a singular direction.

This interpretation is in Stokes' paper 
\cite{stokes1857}. 
Stokes viewed the singular directions as limits of the singular points in the example above, 
as the radius of convergence goes to zero.
Since we are working topologically we can pull these singularities slightly out of the pole, and thus define the tangential punctures
(this is justified by the fact that an equivalence of categories 
can be proved).

The fact that any formal series solution of a linear differential equation is multisummable (and the singular directions are easy to determine) is remarkable.
Multisummation is a morphism of differential algebras 
and so formal solutions sum to actual solutions.

In fact Stokes worked with optimal truncation, not the Borel sum,
so  another puzzle is how to relate them.
For this consider  another familiar  statement:
\beq\begin{matrix}\label{eq: optimal truncation}
\text{\it 
Sometimes the partial sum of a power series 
differs from the full sum}\\
\text{\it by 
less than the modulus of the first term omitted}
\end{matrix}
\eeq
This leads to ``optimal truncation'': 
stopping the sum at the smallest term, and using that to approximate the actual sum.
Stokes applied this to a divergent series and in this way was able to detect the preferred solutions.
This application can be justified
since, in the cases where Stokes applied it,
the statement 
\eqref{eq: optimal truncation} is true provided one replaces ``full sum'' by ``Borel sum'' 
(cf. \cite{watsonbessel} p.219).
Thus the optimal truncation used by Stokes approximates
the Borel sum in much the same way that optimal truncation approximates the usual sum in the case of a convergent series.
The difference in the divergent case is that 
there are singular directions, and the Borel sums on each side of a singular direction are not analytic continuations of each other across the direction.
These are the singular directions detected by Stokes.

\noindent{\bf Acknowledgments.}
The whole philosophy  of {\em deferral of choices}, to set up an intrinsic topological approach, is no doubt due to Deligne---the approaches to Stokes decompositions and Stokes
local systems here are leveraged from his approach to the Stokes filtrations (and their associated graded local systems).  
The Stokes local systems (or {\em wild local systems}) arose from the desire to remove the basepoints in Ramis' wild fundamental group 
\cite{MR91}, as well, of course, 
from the desire to define the local system 
whose monodromy yields the explicit 
presentation given by Birkhoff \cite{birkhoff-1913} (and \cite{JMU81}).
The author lectured on the three main paradigms described here at the Aberdeen ARTIN meeting in March 2018 and the Timi\cb{s}oara GaP meeting in May 2018, and this helped clarify the presentation.
The picture on the title page first appeared in \cite{poster17}.

\renewcommand{\baselinestretch}{1}              %
\normalsize
\bibliographystyle{amsplain}    \label{biby}
\bibliography{../thesis/syr}

\end{document}